\documentclass[12pt]{iopart}
\usepackage{iopams}

\usepackage{graphicx}

\eqnobysec

\begin{document}

\title[Convexification For The Viscosity Solution]{Convexification for the Viscocity Solution for a Coefficient Inverse Problem for the Radiative Transfer Equation}

\author{Michael V. Klibanov$^1$, Jingzhi Li$^2$ and Zhipeng Yang$^3$}

\address{$^1$ Department of Mathematics and Statistics, University of North Carolina at Charlotte, Charlotte, NC, 28223, USA}
\address{$^2$ Department of Mathematics \& National Center for Applied Mathematics Shenzhen \& SUSTech International Center for Mathematics, Southern University of Science and Technology, Shenzhen 518055, P.~R.~China}
\address{$^3$ Department of Mathematics, Southern University of Science and Technology, Shenzhen 518055, P.~R.~China}
\eads{\mailto{mklibanv@uncc.edu}, \mailto{li.jz@sustech.edu.cn}, \mailto{yangzp@sustech.edu.cn}}

\vspace{10pt}
\begin{indented}
\item[]March 2023
\end{indented}

\begin{abstract}
A Coefficient Inverse Problem for the radiative transport equation is
considered. The globally convergent numerical method, the so-called
convexification, is developed. For the first time, the viscosity solution is
considered for a boundary value problem for the resulting system of two
coupled partial differential equations. A Lipschitz stability estimate is
proved for this boundary value problem using a Carleman estimate for the
Laplace operator. Next, the global convergence analysis is provided via that
Carleman estimate. Results of numerical experiments demonstrate a high
computational efficiency of this approach.
\end{abstract}


\submitto{\IP}

%
%

\section{Introduction}

\label{sec:1}

The radiative transfer equation (RTE) is commonly used in diffusive optics
\cite{Heino}. In the case of light propagation, RTE governs scattering and
absorption of photons when they propagate through a diffusive medium, such
as, e.g. turbulent atmosphere and biological medium \cite{Heino}. RTE is
also known as the Boltzmann equation \cite{Chand,Heino}.

One of the specific applications is in seeing through a turbulent
atmosphere. Another attractive application is in optical molecular imaging
(OMI) \cite{W} when some optical markers are attached to specific molecules
to detect faulty genes. In the single positron emission computed tomography
(SPECT) and in photon emission tomography (PET) markers with the high energy
are used. Unlike these, markers of OMI emit a relatively low energy near
infrared light. Reconstruction of the attenuation coefficient using
measurements of the intensity of the emitted light on parts of the human
body should lead to the detection and classification of faulty genes. The
latter might result in a better diagnostics.

Coefficient Inverse Problems (CIPs) for PDEs are both nonlinear and
ill-posed. These two factors cause major challenges in their numerical
solutions. The majority of numerical methods for CIPs is based on a very
popular procedure of the minimization of least squares cost functionals.In
this regard we refer to e.g. \cite{B1,B2,B3,B4,Gonch1,Gonch2,Giorgi,Hassi}. On
the other hand, there is no guarantee of the convexity of these functionals,
and this might lead to the existence of multiple local minima and ravines,
see, e.g. \cite{Scales}.

To avoid the latter, the convexification method was originally proposed in
\cite{KI,Klib97} for the case of CIPs for hyperbolic PDEs. The
convexification is a numerical development of the idea of \cite{BukhKlib},
which was originally proposed only for the proofs of global uniqueness
theorems for multidimensional CIPs and has been explored by many authors
since then, see, e.g. \cite{GY,KP,Lay} as well as \cite{Ksurvey} for a
survey of results as of 2013, \cite{Kpar} for the most recent result, and
the book \cite[Chapters 2,3]{KL}. In \cite{BukhKlib} the method of Carleman
estimates was introduced in the field of CIPs for the first time. Carleman
estimates are also the key to the convexification idea.

In the past few years, members of this research group have been working on
various applications of the convexification to numerical solutions of
various CIPs, see, e.g. \cite{KL} for a summary of results as of 2021. In
particular, in two most recent publications they have solved by
convexification CIPs for two versions of the RTE. \cite{KTR,RRTE}. In both
these works one obtains first a boundary value problem (BVP) for a nonlinear
PDE of the second order, which does not contain the unknown coefficient.
Next, the solution of this BVP is represented via a truncated Fourier series
with $N>2$ terms with respect to a special orthonormal basis. This basis was
originally proposed in \cite{Klib2017}, also, see \cite[section 6.2.3]{KL}.
As a result, one obtains a new BVP for a system of $N$ coupled nonlinear
PDEs of the first order. This BVP is solved via the construction and
subsequent minimization of a weighted cost functional. The main element of
this functional is the presence of a Carleman Weighted Function (CWF) in it.

The CWF is a function, which is used as a weight in the Carleman estimate
for the corresponding PDE operator. The key property of this functional is
that it is strictly convex on a convex bounded set of an arbitrary diameter $%
d>0$ in an appropriate Hilbert space. Since \ a smallness condition is not
imposed on $d$, then it is natural to call that functional \emph{globally}
strictly convex.

The goal of this work is to numerically solve the above mentioned BVP for
that nonlinear PDE via perturbing that PDE operator with the viscosity term $%
-\varepsilon \Delta ,$ where $\varepsilon >0$ is a small number$.$ Then we
obtain a BVP for a system of only two coupled nonlinear PDEs of the second
order, as opposed to the above case of $N>2$ PDEs . First, we prove the
Lipschitz stability estimate for this BVP using a Carleman estimate for the
Laplace operator. Next, we construct a globally strictly convex weighted
cost functional for the latter BVP. The key element of this functional is
the presence of the Carleman Weight Function in it. This function is
involved as the weight in the above mentioned Carleman estimate. The proof
of the global strict convexity of that functional is the main part of our
global convergence analysis. This analysis ends up with the proof of the
global convergence of the gradient descent method of the minimization of
that functional. Finally, we conduct exhaustive numerical studies, which
demonstrate a good performance of our technique. Previously the viscosity
solution of the Hamilton-Jacobi equation was numerically computed in \cite%
{KHJ} by a different version of the convexification method.

Publications \cite{KTR,RRTE} represent the first numerical solutions of CIPs
for both RTE \cite{KTR} and its Riemannian version \cite{RRTE}. Previous
numerical results were obtained only for the inverse source problems for RTE
in \cite{HT1,HT2,HT3,Smirnov} and for RRTE in \cite{Le}. However, inverse
source problems are linear ones, unlike CIPs. There are also some uniqueness
and stability results for CIPs for RTE. Since this paper is not concerned
with such results, then we refer here only to a limited number of
publications on the latter topic \cite{Bal,GY,KP,Lay}.

In section 2 we state both forward and inverse problems. To solve the
inverse problem, we obtain in in section 3 a boundary value problem for a $%
2\times 2$ system of nonlinear PDEs with viscosity terms. We prove the
Lipschitz stability estimate for this boundary value problem in section 4.
In section 5 we construct the convexified Tikhonov-like functional and
conduct the convergence analysis for it. Numerical experiments are presented
in section 6.

\section{Statements of Forward and Inverse Problems}

\label{sec:2}

For $n\geq 1,$ points in $\mathbb{R}^{n}$ are denoted below as $%
x=(x_{1},x_{2},...,x_{n-1},z)\in \mathbb{R}^{n}.$ Let numbers $B,a,b,d>0$,
where
\begin{equation}
0<a<b.  \label{2.0}
\end{equation}%
Let $\Omega \subset \mathbb{R}^{n}$ be a rectangular prism with its boundary
$\partial \Omega =\partial _{1}\Omega \cup \partial _{2}\Omega \cup \partial
_{3}\Omega ,$ where
\begin{eqnarray}
& \Omega =\{x:-B<x_{1},...,x_{n-1}<B,a<z<b\},  \label{2.1} \\
& \hspace{0.2cm}\partial _{1}\Omega =\left\{
x:-B<x_{1},...,x_{n-1}<B,z=a\right\} ,  \label{2.2} \\
& \hspace{0.2cm}\partial _{2}\Omega =\left\{ x\mathbf{:}
-B<x_{1},...,x_{n-1}<B,z=b\right\} ,  \label{2.3} \\
& \hspace{0.1cm}\partial _{3}\Omega =\left\{ x_{i}=\pm B,z\in \left(
a,b\right) ,i=1,...,n-1\right\} .  \label{2.4}
\end{eqnarray}%
Let $\Phi _{d}$ be the line, where the external sources are located,
\begin{equation}
\Phi _{d}=\{x_{\alpha }=(\alpha ,0,...,0):\alpha \in \lbrack -d,d]\}.
\label{2.40}
\end{equation}%
Hence, $\Phi _{d}$ is a part of the $x_{1}-$axis. By (\ref{2.0}), (\ref{2.1}%
) and (\ref{2.40})
\begin{equation}
\Phi _{d}\cap \overline{\Omega }=\varnothing .  \label{2.400}
\end{equation}

Let $\sigma >0$ be a sufficiently small number. To simplify the
presentation, avoid working with singularities. Hence, we consider the
following function instead of $\delta \left( x\right) $:%
\begin{equation}
f_{\sigma }\left( x\right) =C_{\sigma}\left\{
\begin{array}{c}
\exp \left( \frac{\left\vert x\right\vert ^{2}}{\sigma ^{2}-\left\vert
x\right\vert ^{2}}\right) ,\left\vert x\right\vert <\sigma , \\
0,\left\vert x\right\vert \geq \sigma%
\end{array}
\right. .  \label{2.5}
\end{equation}%
The constant $C_{\sigma }$ is chosen here such that
\begin{equation}
C_{\sigma }\int\limits_{\left\vert x\right\vert <\sigma }\exp \left( \frac{
\left\vert x\right\vert ^{2}}{\sigma ^{2}-\left\vert x\right\vert ^{2}}
\right) dx=1.  \label{2.50}
\end{equation}%
Hence, the function $f\left( x-x_{\alpha }\right) =f\left( x_{1}-\alpha
,x_{2},...,x_{n},z\right) \in C^{\infty }\left( \mathbb{R}^{n}\right) $ can
be considered as the source function for the source $x_{\alpha }\in \Phi
_{d}.$ Let the number $\sigma $ be so small that
\begin{equation}
f\left( x-x_{\alpha }\right) =0,\forall x\in \overline{\Omega },\forall
x_{\alpha }\in \Phi _{d}.  \label{2.6}
\end{equation}

Let
\begin{equation}
\overline{B}=\max \left( B,d\right) .  \label{2.60}
\end{equation}
Introduce the domain $P\subset \mathbb{R}^{n}$ as well as two of its
subdomains $P_{a}^{+}$,$P_{a}^{-}\subset P,$%
\begin{equation}
\left\{
\begin{array}{c}
P=\left\{ x:\left\vert x_{1}\right\vert ,...,\left\vert x_{n-1}\right\vert <
\overline{B},z\in \left( 0,b\right) \right\} , \\
P_{a}^{+}=P\cap \left\{ z>a\right\} , \\
P_{a}^{-}=P\diagdown P_{a}^{+}.%
\end{array}
\right.  \label{2.07}
\end{equation}%
By (\ref{2.1}), (\ref{2.60}) and (\ref{2.07})
\begin{equation}
\Omega \subseteq P_{a}^{+}.  \label{2.007}
\end{equation}%
Below%
\begin{equation}
x\in P, \ \alpha \in \left( -d,d\right) .  \label{2.170}
\end{equation}

For two arbitrary points $x,y\in \mathbb{R}^{n}$ denote $L\left( x,y\right) $
the line segment connecting these two points and let $ds$ be the element of
the euclidean length on $L\left( x,y\right) .$ Let $\nu (x\mathbf{,}\alpha )$
be the unit vector, which is parallel to $L\left( x,x_{\alpha }\right) ,$
\begin{equation}
\nu (x\mathbf{,}\alpha )=\frac{x-x_{\alpha }}{\left\vert x-x_{\alpha
}\right\vert },\ x\neq x_{\alpha }.  \label{2.130}
\end{equation}

Let $u(x,\alpha )$ be the steady-state radiance at the point $x$ generated
by the source function $f\left( x-x_{\alpha }\right) .$ The function $%
u(x,\alpha )$ satisfies the stationary RTE \cite{Heino}:
\begin{equation}
\left.
\begin{array}{c}
\nu (x\mathbf{,}\alpha )\cdot \nabla _{x}u(x,\alpha )+a\left( x\right)
u(x,\alpha ) \\
=\mu _{s}(x)\int\limits_{\Phi _{d}}G(x,\alpha ,\beta )u(x,\beta )d\beta
+f\left( x-x_{\alpha }\right) ,\ x\in P,x_{\alpha }\in \Phi _{d},%
\end{array}
\right.  \label{2.7}
\end{equation}
where \textquotedblleft $\cdot $" denotes the scalar product in $\mathbb{R}%
^{n}$, see (\ref{2.170}). The kernel $G(x,\alpha ,\beta )$ of the integral
operator in (\ref{2.7}) satisfies \cite{Heino,Smirnov}:
\begin{eqnarray}
&G(x,\alpha ,\beta )\geq 0,\ x\mathbf{\in }\overline{P};\ \alpha
,\beta \in \left[ -d,d\right] ,\   \label{2.8} \\
&\hspace{0.5 cm} G(x,\alpha ,\beta )\in C^{1}\left( \overline{P}\times \left[
-d,d\right] ^{2}\right) .  \label{2.9}
\end{eqnarray}
In (\ref{2.7}),%
\begin{equation}
a\left( x\right) =\mu _{a}\left( x\right) +\mu _{s}(x),  \label{2.10}
\end{equation}%
where $\mu _{a}\left( x\right) $ and $\mu _{s}(x)$ are the absorption and
scattering coefficients respectively. The function $a\left( x\right) $ is
the attenuation coefficient. We assume that
\begin{eqnarray}
& \mu _{a}\left( x\right) ,\mu _{s}(x)\geq 0,\mu _{a}\left( x\right) =\mu
_{s}(x)=0,\  x\in \mathbb{\ }P\setminus \Omega ,  \label{2.11} \\
& \hspace{2cm}\mu _{a}\left( x\right) ,\mu _{s}(x)\in C\left( \overline{P}
\right) .  \label{2.12}
\end{eqnarray}

\textbf{Forward Problem.} \emph{Let conditions (\ref{2.0})-(\ref{2.12})
hold. Find the function }$u(x,\alpha )$ $\in C^{1}\Big(\overline{P}\times %
\left[ -d,d\right] \Big)$\emph{\ satisfying equation (\ref{2.7}) and the
initial condition}%
\begin{equation}
u(x_{\alpha },\alpha )=0\mbox{ for }x_{\alpha }\in \Phi _{d}.  \label{2.14}
\end{equation}

Denote
\begin{equation*}
c\left( x\right) =\exp \left( \int_{L(x,x_{\alpha })}a(x\left( s\right)
)ds\right) .
\end{equation*}%
The following existence and uniqueness theorem for the Forward Problem was
proven in \cite{KTR}, and a similar theorem was proven in \cite{RRTE} for
the Riemannian analog of RTE:

\textbf{Theorem 2.1 }\cite{KTR}. \emph{Assume that conditions (\ref{2.0})-( %
\ref{2.130}) and (\ref{2.8})-(\ref{2.12}) hold. Then} \emph{there exists
unique solution }$u\left( x,\alpha \right) \in C^{1}\left( \overline{P}%
\times \left[ -d,d\right] \right) $\emph{\ of equation (\ref{2.7}) with the
initial condition (\ref{2.14}), and the function }$u\left( x,\alpha \right) $%
\emph{\ has the following form for }$x\in P,\alpha \in \left( -d,d\right) :$%
\begin{eqnarray}
& \hspace{3cm}u(x\mathbf{,}\alpha )=u_{0}(x\mathbf{,}\alpha )+  \label{2.025}
\\
& \hspace{-1cm} +\frac{1}{c(x\mathbf{,}\alpha )}\int\limits_{L(x,x_{\alpha })}c(x\left(
s\right) ,\alpha )\mu _{s}(x\left( s\right) )\left( \int\limits_{\Phi
_{d}}G(x\left( s\right) ,\alpha ,\beta )u(x\left( s\right) ,\beta )d\beta
\right) ds,  \nonumber \\
& \hspace{2cm} u_{0}(x\mathbf{,}\alpha )=\frac{1}{c\left( x,\alpha \right) }
\int\limits_{L(x,x_{\alpha })}f\left( x\left( s\right) -x_{\alpha }\right)
ds.  \label{2.026}
\end{eqnarray}%
\emph{\ Furthermore, the following inequality holds: }%
\begin{eqnarray}
u\left( x,\alpha \right) \geq m>0\mbox{ for }\left( x,\alpha \right) \in
\overline{P}_{a}^{+}\times \left[ -d,d\right] , \ m=\min_{\overline{P}%
_{a}^{+}\times \left[ -d,d\right] }u_{0}(x \mathbf{,}\alpha ).  \label{2.16}
\end{eqnarray}

\textbf{Coefficient Inverse Problem (CIP).} \emph{Let conditions (\ref{2.0}
)-(\ref{2.12}) hold. Let the function }$u\left( x,\alpha \right) \in
C^{1}\left( \overline{G}\times \left[ -d,d\right] \right) $\emph{\ be the
solution of the Forward Problem as in Theorem 2.1. Assume that the
attenuation coefficient }$a\left( x\right) $\emph{\ in (\ref{2.7}) is
unknown. Assume that the function }$g\left( x,\alpha \right) $\emph{\ is
known,}
\begin{equation}
g\left( x,\alpha \right) =u\left( x,\alpha \right) ,\forall x\in \partial
\Omega ,\forall \alpha \in \left( -d,d\right) .  \label{2.15}
\end{equation}%
\emph{Find the function }$a\left( x\right) .$\emph{\ }

\section{Boundary Value Problem for a $2\times 2$ System of PDEs With
Viscosity Terms}

\label{sec:3}

\subsection{Preliminaries}

\label{sec:3 copy(1)}

It follows from (\ref{2.007}), (\ref{2.170}) and (\ref{2.16}) that we can
introduce a new function $w(x,\alpha )$,
\begin{equation}
w\left( x,\alpha \right) =\ln u\left( x,\alpha \right) ,\ \left(
x,\alpha \right) \in \overline{\Omega }\times \left[ -d,d\right] .
\label{3.1}
\end{equation}%
Substituting (\ref{3.1}) in (\ref{2.7}) and (\ref{2.15}), we obtain
\begin{eqnarray}
& \hspace{2.5cm}\nu (x\mathbf{,}\alpha )\cdot \nabla _{x}w(x,\alpha
)+a\left( x\right)  \label{3.2} \\
& =e^{-w\left( x,\alpha \right) }\mu _{s}(x)\int\limits_{\Phi
_{d}}G(x,\alpha ,\beta )e^{w\left( x,\beta \right) }d\beta ,x\in \Omega
,\alpha \in \left( -d,d\right) ,  \nonumber \\
& \hspace{2cm}w(x,\alpha )\mid _{\partial \Omega }=g_{1}\left( x,\alpha
\right) =\ln g\left( x,\alpha \right) .  \label{3.3}
\end{eqnarray}

Differentiate both sides of (\ref{3.2}) with respect to $\alpha $ and use $%
\partial _{\alpha }a\left( x\right) \equiv 0.$ We obtain an integral
differential equation with the derivatives up to the second order,
\begin{eqnarray}
&\hspace{1 cm} \nu (x\mathbf{,}\alpha )\cdot \nabla _{x}w_{\alpha }(x,\alpha
)+\partial _{\alpha }\nu (x\mathbf{,}\alpha )\cdot \nabla _{x}w(x,\alpha )+
\\
&\hspace{0.8 cm} +e^{-w\left( x,\alpha \right) }w_{\alpha }\left( x,\alpha
\right) \mu _{s}(x)\int\limits_{\Phi _{d}}G(x,\alpha ,\beta )e^{w\left(
x,\beta \right) }d\beta -  \label{3.5} \\
&-e^{-w\left( x,\alpha \right) }\mu _{s}(x)\int\limits_{\Phi _{d}}\partial
_{\alpha }G(x,\alpha ,\beta )e^{w\left( x,\beta \right) }d\beta =0,x\in
\Omega ,\alpha \in \left( -d,d\right) .
\end{eqnarray}

By (\ref{3.3}) the boundary condition for $w_{\alpha }(x,\alpha )$ is
\begin{equation}
w_{\alpha }(x,\alpha )\mid _{\partial \Omega }=g_{2}\left( x,\alpha \right)
= \frac{g_{\alpha }\left( x,\alpha \right) }{g\left( x,\alpha \right) }.
\label{3.6}
\end{equation}

To work with the viscosity solution, we need to figure out the following
boundary conditions, see (\ref{2.3}):%
\begin{eqnarray}
& \hspace{0.2 cm} \partial _{z}w(x,\alpha )\mid _{\partial _{2}\Omega
}=\partial _{z}w\left( x_{1},...,x_{n-1},b,\alpha \right) =g_{3}(x,\alpha ),
\label{3.60} \\
& \partial _{z}w_{\alpha }(x,\alpha )\mid _{\partial _{2}\Omega }=\partial
_{z}w_{\alpha }\left( x_{1},...,x_{n-1},b,\alpha \right) =g_{4}(x,\alpha ).
\label{3.61}
\end{eqnarray}
It follows from (\ref{2.40}) and (\ref{2.130}) that $\nu (x\mathbf{,}\alpha
)=( \nu _{1}(x\mathbf{,}\alpha )$, $\nu _{2}(x\mathbf{,}\alpha )$, $\cdots$,
$\nu_{n-1}(x\mathbf{,}\alpha )$, $\nu _{n}(x\mathbf{,}\alpha )) $, where
\begin{eqnarray}
&\hspace{1 cm} \nu _{1}(x\mathbf{,}\alpha )=\frac{\left( x_{1}-\alpha
\right) }{\sqrt{ \left( x_{1}-\alpha \right)
^{2}+x_{2}^{2}+...+x_{n-1}^{2}+z^{2}}},  \label{3.62} \\
& \nu _{k}(x\mathbf{,}\alpha )=\frac{x_{k}}{\sqrt{\left( x_{1}-\alpha
\right) ^{2}+x_{2}^{2}+...+x_{n-1}^{2}+z^{2}}},k=2,...,n-1,  \label{3.63} \\
&\hspace{1 cm} \nu _{n}(x\mathbf{,}\alpha )=\frac{z}{\sqrt{\left(
x_{1}-\alpha \right) ^{2}+x_{2}^{2}+...+x_{n-1}^{2}+z^{2}}}.  \label{3.64}
\end{eqnarray}
Hence,
\begin{equation}
\nu _{n}(x\mathbf{,}\alpha )\mid _{\partial _{2}\Omega }=\frac{b}{\sqrt{
\left( x_{1}-\alpha \right) ^{2}+x_{2}^{2}+...+x_{n-1}^{2}+b^{2}}}.
\label{3.65}
\end{equation}%
By (\ref{2.10}) and (\ref{2.11}) $a\left( x\right) =0$ for $x\in \partial
_{2}\Omega .$ Hence, (\ref{3.2}), (\ref{3.3}) and (\ref{3.60})-(\ref{3.65})
imply:%
\begin{eqnarray}
&\hspace{2.8 cm} \partial _{z}w(x,\alpha )\mid _{\partial _{2}\Omega
}=g_{3}(x,\alpha )= \\
&\hspace{-1 cm} -\frac{1}{\nu _{n}(x\mathbf{,}\alpha )}\left[ \sum\limits_{k=1}^{n-1}\nu
_{k}(x\mathbf{,}\alpha )w_{x_{k}}(x\mathbf{,}\alpha )+\frac{\mu _{s}(x)}{
g\left( x,\alpha \right) }\int\limits_{\Phi _{d}}G(x,\alpha ,\beta )g\left(
x,\beta \right) d\beta \right] ,  \label{3.66} \\
&\hspace{3.2 cm} g_{4}(x,\alpha )=\partial _{\alpha }g_{3}(x,\alpha ),
\label{3.67} \\
&\hspace{3.2 cm} x\in \partial _{2}\Omega ,\alpha \in \left( -d,d\right) .
\label{3.68}
\end{eqnarray}

\subsection{Viscosity solution}

\label{sec:2.2}

Denote%
\begin{equation}
p\left( x,\alpha \right) =w\left( x,\alpha \right) ,\ q\left( x,\alpha
\right) =w_{\alpha }\left( x,\alpha \right) .  \label{3.7}
\end{equation}%
Based on (\ref{3.3}), (\ref{3.5}), (\ref{3.6}), (\ref{3.66})-(\ref{3.68}),
consider the following BVP for the system of viscosity equations with a
small parameter $\varepsilon >0$:%
\begin{equation}
\hspace{-1 cm}
\left.
\begin{array}{c}
L_{1}\left( p_{\varepsilon },q_{\varepsilon }\right) =-\varepsilon \Delta
p_{\varepsilon }+\nu (x\mathbf{,}\alpha )\cdot \nabla _{x}q_{\varepsilon
}(x,\alpha )+\partial _{\alpha }\nu (x\mathbf{,}\alpha )\cdot \nabla
_{x}p_{\varepsilon }(x,\alpha )+ \\
+e^{-p_{\varepsilon }\left( x,\alpha \right) }q_{\varepsilon }(x,\alpha )\mu
_{s}(x)\int\limits_{\Phi _{d}}G(x,\alpha ,\beta )e^{p_{\varepsilon }\left(
x,\beta \right) }d\beta - \\
-e^{-p_{\varepsilon }\left( x,\alpha \right) }\mu _{s}(x)\int\limits_{\Phi
_{d}}\partial _{\alpha }G(x,\alpha ,\beta )e^{p_{\varepsilon }\left( x,\beta
\right) }d\beta =0,\ x\in \Omega ,\alpha \in \left( -d,d\right) ,%
\end{array}
\right.  \label{3.8}
\end{equation}%
And also%
\begin{eqnarray}
&\hspace{-1 cm} \left.
\begin{array}{c}
L_{2}\left( p_{\varepsilon },q_{\varepsilon }\right) =-\varepsilon \Delta
q_{\varepsilon }+\nu (x\mathbf{,}\alpha )\cdot \nabla _{x}q_{\varepsilon
}(x,\alpha )+\partial _{\alpha }\nu (x\mathbf{,}\alpha )\cdot \nabla
_{x}p_{\varepsilon }(x,\alpha )+ \\
+e^{-p_{\varepsilon }\left( x,\alpha \right) }q_{\varepsilon }(x,\alpha )\mu
_{s}(x)\int\limits_{\Phi _{d}}G(x,\alpha ,\beta )e^{p_{\varepsilon }\left(
x,\beta \right) }d\beta - \\
-e^{-p_{\varepsilon }\left( x,\alpha \right) }\mu _{s}(x)\int\limits_{\Phi
_{d}}\partial _{\alpha }G(x,\alpha ,\beta )e^{p_{\varepsilon }\left( x,\beta
\right) }d\beta ,\ x\in \Omega ,\alpha \in \left( -d,d\right) .%
\end{array}
\right.  \label{3.9} \\
&\hspace{2 cm} p_{\varepsilon }\mid _{\partial \Omega }=g_{1}\left( x,\alpha
\right) ,\partial _{z}p_{\varepsilon }\mid _{\partial _{2}\Omega
}=g_{3}\left( x,\alpha \right) ,  \label{3.10} \\
&\hspace{2 cm} q_{\varepsilon }\mid _{\partial \Omega }=g_{2}\left( x,\alpha
\right) ,\partial _{z}q_{\varepsilon }\mid _{\partial _{2}\Omega
}=g_{4}\left( x,\alpha \right) .  \label{3.11}
\end{eqnarray}

Therefore, we have obtained the BVP (\ref{3.8})-(\ref{3.11}) with respect to
the pair of functions $\left( p_{\varepsilon },q_{\varepsilon }\right)
\left( x,\alpha \right) $. By (\ref{3.10}) and (\ref{3.11}) this BVP\ has an
overdetermination in the Neumann boundary conditions at $z=b.$ We focus
below on the solution of this BVP. Suppose that we have computed a solution $%
\left( p_{\varepsilon ,\mbox{comp}},q_{\varepsilon ,\mbox{comp}}\right)
\left( x,\alpha \right) $ of BVP (\ref{3.8})-(\ref{3.11}). Then we use (\ref%
{3.2}) and (\ref{3.7}) to compute the target coefficient $a\left( x\right) ,$
\begin{equation}
\hspace{-1 cm}
\left.
\begin{array}{c}
a_{\mbox{comp}}\left( x\right) =-\frac{1}{2d}\int\limits_{-d}^{d}\nu (x
\mathbf{,}\alpha )\cdot \nabla _{x}p_{\varepsilon ,\mbox{comp}}(x,\alpha
)d\alpha + \\
+\frac{1}{2d}\int\limits_{-d}^{d}e^{-p_{\varepsilon ,\mbox{comp}}\left(
x,\alpha \right) }\left( \mu _{s}(x)\int\limits_{\Phi _{d}}G(x,\alpha ,\beta
)e^{p_{\varepsilon ,\mbox{comp}}\left( x,\beta \right) }d\beta \right)
d\alpha ,\ x\in \Omega .%
\end{array}
\right.  \label{3.12}
\end{equation}

\textbf{Remark 3.1.} \emph{As soon as we got BVP (\ref{3.8})-(\ref{3.11})
for the pair of functions }$\left( p_{\varepsilon },q_{\varepsilon }\right)
, $\emph{\ we do not require anymore that }%
\begin{equation}
q_{\varepsilon }=\partial _{\alpha }p_{\varepsilon }  \label{3.120}
\end{equation}%
\emph{\ as in (\ref{3.7}). In other words, we solve this BVP for a slightly
broader class of vector functions }$\left( p_{\varepsilon },q_{\varepsilon
}\right) .$\emph{\ Nevertheless, it follows from the uniqueness claim of
Theorem 4.1 that the solution of this \textquotedblleft broader" BVP, if it
exists, is still such that (\ref{3.120}) holds. }

\section{Lipschitz Stability Estimate for BVP (\protect\ref{3.8})-(\protect
\ref{3.11})}

\label{sec:4}

Introduce the space $H_{2}^{1}\left( \Omega \right) \times L_{2,2}\left(
-d,d\right) $ as
\begin{equation}
\left.
\begin{array}{c}
H_{2}^{1}\left( \Omega \right) \times L_{2,2}\left( -d,d\right) = \\
= \left\{
\begin{array}{c}
\left( p,q\right) :\left\Vert \left( p,q\right) \right\Vert
_{H_{2}^{1}\left( \Omega \right) \times L_{2,2}\left( -d,d\right) }^{2}= \\
= \int\limits_{-d}^{d}\left( \left\Vert p\left( x,\alpha \right) \right\Vert
_{H^{1}\left( \Omega \right) }^{2}+\left\Vert q\left( x,\alpha \right)
\right\Vert _{H^{1}\left( \Omega \right) }^{2}\right) d\alpha <\infty .%
\end{array}
\right\}%
\end{array}
\right.  \label{5.002}
\end{equation}

\textbf{Theorem 4.1} (Lipschitz stability and uniqueness). \emph{Assume that
conditions (\ref{2.0})-(\ref{2.400}) and (\ref{2.9})-(\ref{2.12}) hold.
Suppose that there exists two pairs of functions }

$\left( p_{\varepsilon ,1},q_{\varepsilon ,1}\right) \left( x,\alpha \right)
,\left( p_{\varepsilon ,2},q_{\varepsilon ,2}\right) \left( x,\alpha \right)
\in H^{2}\left( \Omega \right) \times C\left[ -d,d\right] $\emph{\
satisfying equations (\ref{3.8}) and (\ref{3.9}) and such that }%
\begin{equation}
\hspace{-1 cm}
p_{\varepsilon ,1}\left( x,\alpha \right) =p_{\varepsilon ,2}\left( x,\alpha
\right) , \ q_{\varepsilon ,1}\left( x,\alpha \right) =q_{\varepsilon
,2}\left( x,\alpha \right), \ x\in \partial \Omega \diagdown \partial
_{2}\Omega , \ \alpha \in \left( -d,d\right) .
\label{4.2}
\end{equation}%
\emph{Let }%
\begin{equation*}
M=\max \left\{ \max_{i=1,2}\left\Vert p_{\varepsilon ,i}\right\Vert
_{C\left( \overline{\Omega }\times \left[ -d,d\right] \right)
},\max_{i=1,2}\left( \left\Vert q_{\varepsilon ,i}\right\Vert _{C\left(
\overline{\Omega }\times \left[ -d,d\right] \right) }\right) \right\} .
\end{equation*}%
\emph{\ Then there exists a constant}%
\begin{equation}
C_{1}=C_{1}\left( \Omega ,d,\varepsilon ,\left\Vert G\right\Vert _{C\left(
\overline{\Omega }\right) \times C^{1}\left[ -d,d\right] ^{2}},\left\Vert
\mu _{a}\right\Vert _{C\left( \overline{\Omega }\right) },\left\Vert \mu
_{s}\right\Vert _{C\left( \overline{\Omega }\right) },M\right) >0
\label{4.02}
\end{equation}%
\emph{depending only on listed parameters such that the following Lipschitz
stability estimates hold:}
\begin{equation}
\hspace{-1 cm}
\left.
\begin{array}{c}
\left\Vert \left( p_{\varepsilon ,1}-p_{\varepsilon ,2},q_{\varepsilon
,1}-q_{\varepsilon ,2}\right) \right\Vert _{H_{2}^{1}\left( \Omega \right)
\times L_{2,2}\left( -d,d\right) }\leq \\
\leq C_{1}\left( \left\Vert p_{\varepsilon ,1}-p_{\varepsilon ,2}\right\Vert
_{H^{1}\left( \partial _{2}\Omega \right) \times L_{2}\left( -d,d\right)
}+\left\Vert \partial _{z}p_{\varepsilon ,1}-\partial _{z}p_{\varepsilon
,2}\right\Vert _{L_{2}\left( \partial _{2}\Omega \right) \times L_{2}\left(
-d,d\right) }\right) + \\
+C_{1}\left( \left\Vert q_{\varepsilon ,1}-q_{\varepsilon ,2}\right\Vert
_{H^{1}\left( \partial _{2}\Omega \right) \times L_{2}\left( -d,d\right)
}+\left\Vert \partial _{z}q_{\varepsilon ,1}-\partial _{z}q_{\varepsilon
,2}\right\Vert _{L_{2}\left( \partial _{2}\Omega \right) \times L_{2}\left(
-d,d\right) }\right) .%
\end{array}
\right.  \label{4.3}
\end{equation}
\emph{Let }$a_{\varepsilon ,1}\left( x\right) $\emph{\ and }$a_{\varepsilon
,2}\left( x\right) $ \emph{\ functions computed via the right hand side of ( %
\ref{3.12}), in which }$p_{\varepsilon ,\mbox{comp}}(x,\alpha )$\emph{\ is
replaced with }$p_{\varepsilon ,1}(x,\alpha )$\emph{\ and }$p_{\varepsilon
,2}(x,\alpha )$\emph{\ respectively. Then the following Lipschitz stability
estimate holds:}
\begin{equation}
\hspace{-1 cm}
\left.
\begin{array}{c}
\left\Vert a_{\varepsilon ,1}-a_{\varepsilon ,2}\right\Vert _{L_{2}\left(
\Omega \right) }\leq \\
\leq C_{1}\left( \left\Vert p_{\varepsilon ,1}-p_{\varepsilon ,2}\right\Vert
_{H^{1}\left( \partial _{2}\Omega \right) \times L_{2}\left( -d,d\right)
}+\left\Vert \partial _{z}p_{\varepsilon ,1}-\partial _{z}p_{\varepsilon
,2}\right\Vert _{L_{2}\left( \partial _{2}\Omega \right) \times L_{2}\left(
-d,d\right) }\right) + \\
+C_{1}\left( \left\Vert q_{\varepsilon ,1}-q_{\varepsilon ,2}\right\Vert
_{H^{1}\left( \partial _{2}\Omega \right) \times L_{2}\left( -d,d\right)
}+\left\Vert \partial _{z}q_{\varepsilon ,1}-\partial _{z}q_{\varepsilon
,2}\right\Vert _{L_{2}\left( \partial _{2}\Omega \right) \times L_{2}\left(
-d,d\right) }\right) .%
\end{array}
\right.  \label{4.4}
\end{equation}
\emph{In particular, suppose that, in addition to (\ref{4.2}) }%
\begin{equation*}
\hspace{-1.5 cm}
p_{\varepsilon ,1}=p_{\varepsilon ,2},q_{\varepsilon ,1}=q_{\varepsilon
,2},\partial _{z}p_{\varepsilon ,1}=\partial _{z}p_{\varepsilon ,2},\partial
_{z}q_{\varepsilon ,1}=\partial _{z}q_{\varepsilon ,2}\mbox{ for }\left(
x,\alpha \right) \in \partial _{2}\Omega \times \left( -d,d\right) .
\end{equation*}%
\emph{Then (\ref{4.3}) and (\ref{4.4}) imply that }$\emph{for\ }x\in \Omega
,\alpha \in \left( -d,d\right) $\emph{\ }%
\begin{equation}
\hspace{-1 cm}
p_{\varepsilon ,1}\left( x,\alpha \right) \equiv p_{\varepsilon ,2}\left(
x,\alpha \right) ,q_{\varepsilon ,1}\left( x,\alpha \right) \equiv
q_{\varepsilon ,2}\left( x,\alpha \right) \emph{\ and\ }a_{\varepsilon
,1}\left( x\right) \equiv a_{\varepsilon ,2}\left( x\right) \emph{.}
\label{1}
\end{equation}

Below $C_{1}>0$ denotes different constants depending on parameters listed
in (\ref{4.02}). Uniqueness of the BVP (\ref{3.8})-(\ref{3.11}) obviously
follows from (\ref{1}).

\subsection{Carleman estimate}

\label{sec:4.1}

Prior the proof of Theorem 4.1, we need to prove a Carleman estimate for the
operator $\Delta .$ \ Denote%
\begin{equation}
H_{0}^{2}\left( \Omega \right) =\left\{ u\in H^{2}\left( \Omega \right)
:u\mid _{\partial \Omega \diagdown \partial _{2}\Omega }=0\right\} .
\label{4.05}
\end{equation}

\textbf{Theorem 4.2. }\emph{Assume that conditions}\textbf{\ }\emph{(\ref%
{2.1})-(\ref{2.4}) hold. There exists a constant }$C=C\left( \Omega \right)
>0$\emph{\ and a sufficiently large number }$\lambda _{0}=\lambda _{0}\left(
\Omega \right) >1$\emph{,} \emph{both depending only on the domain }$\Omega
, $ \emph{such that the following Carleman estimate holds for all \ and for
all }$\lambda \geq \lambda _{0}:$
\begin{equation}
\left.
\begin{array}{c}
\int\limits_{\Omega }\left( \Delta u\right) ^{2}e^{2\lambda z^{2}}dx\geq
C\int\limits_{\Omega }\left( \lambda \left( \nabla u\right) ^{2}+\lambda
^{3}u^{2}\right) e^{2\lambda z^{2}}dx- \\
-C\lambda ^{3}\left( \left\Vert u\right\Vert _{H^{1}\left( \partial
_{2}\Omega \right) }^{2}+\left\Vert \partial _{z}u\right\Vert _{L_{2}\left(
\partial _{2}\Omega \right) }^{2}\right) e^{2\lambda b^{2}},\forall \lambda
\geq \lambda _{0},\forall u\in H_{0}^{2}\left( \Omega \right) .%
\end{array}
\right.  \label{4.5}
\end{equation}

\textbf{Proof.} Everywhere below $C=C\left( \Omega \right) >0$ denotes
different constants depending only on the domain $\Omega .$ We assume first
that
\begin{equation}
u\in C^{2}\left( \overline{\Omega }\right) \cap H_{0}^{2}\left( \Omega
\right) .  \label{4.500}
\end{equation}%
Introduce a new function $v\left( x\right) ,$%
\begin{equation}
v\left( x\right) =u\left( x\right) e^{\lambda z^{2}}.  \label{4.6}
\end{equation}%
Hence,%
\begin{equation}
\hspace{-1 cm}
\left.
\begin{array}{c}
u=ve^{-\lambda z^{2}}, \quad u_{z}=\left( v_{z}-2\lambda zv\right)
e^{-\lambda z^{2}}, \\
u_{zz}=\left( v_{zz}-4\lambda zv_{z}+4\lambda ^{2}\left( z^{2}-z/\left(
2\lambda \right) \right) v\right) e^{-\lambda z^{2}}, \quad
u_{x_{i}x_{i}}=v_{x_{i}x_{i}}e^{-\lambda z^{2}}.%
\end{array}
\right.  \label{4.7}
\end{equation}%
Hence,
\begin{equation}
\hspace{-1 cm}
\left.
\begin{array}{c}
\left( \Delta u\right) ^{2}e^{2\lambda z^{2}}=\left[ \left(
\sum\limits_{i=1}^{n-1}v_{x_{i}x_{i}}+v_{zz}+4\lambda ^{2}\left(
z^{2}-z/\left( 2\lambda \right) \right) v\right) -4\lambda zv_{z}\right]
^{2}\geq \\
\geq -8\lambda zv_{z}\left(
\sum\limits_{i=1}^{n-1}v_{x_{i}x_{i}}+v_{zz}+4\lambda ^{2}\left(
z^{2}-z/\left( 2\lambda \right) \right) v\right) .%
\end{array}
\right.  \label{4.8}
\end{equation}
\textbf{Step 1.} Estimate from the below the following term in the second
line of (\ref{4.8}):
\begin{equation*}
\left.
\begin{array}{c}
-8\lambda zv_{z}\left( \sum\limits_{i=1}^{n-1}v_{x_{i}x_{i}}+v_{zz}\right) =
\\
=\sum\limits_{i=1}^{n-1}\left( -8\lambda zv_{z}v_{x_{i}}\right)
_{x_{i}}+\sum\limits_{i=1}^{n-1}\left( 8\lambda zv_{zx_{i}}v_{x_{i}}\right)
+\left( -4\lambda zv_{z}^{2}\right) _{z}+4\lambda v_{z}^{2}\geq \\
\geq -4\lambda \sum\limits_{i=1}^{n-1}v_{x_{i}}^{2}+\left( -4\lambda
zv_{z}^{2}+4\lambda z\sum\limits_{i=1}^{n-1}v_{x_{i}}^{2}\right)
_{z}+\sum\limits_{i=1}^{n-1}\left( -8\lambda zv_{z}v_{x_{i}}\right) _{x_{i}}.%
\end{array}
\right.
\end{equation*}
Thus, moving from the function $v$ to the function $u$ via (\ref{4.6}), we
obtain
\begin{equation}
\left.
\begin{array}{c}
-8\lambda zv_{z}\left( \sum\limits_{i=1}^{n-1}v_{x_{i}x_{i}}+v_{zz}\right)
\geq -4\lambda \sum\limits_{i=1}^{n-1}u_{x_{i}}^{2}e^{2\lambda z^{2}}+ \\
+\left( -4\lambda z\left( \left( u_{z}+2\lambda zu\right) e^{\lambda
z^{2}}\right) _{z}^{2}+4\lambda
z\sum\limits_{i=1}^{n-1}u_{x_{i}}^{2}e^{2\lambda z^{2}}\right) _{z}+ \\
+\sum\limits_{i=1}^{n-1}\left( -8\lambda z\left( u_{z}+2\lambda zu\right)
u_{x_{i}}e^{2\lambda z^{2}}\right) _{x_{i}}.%
\end{array}
\right.  \label{4.9}
\end{equation}

\textbf{Step 2.} Estimate from the below the following term in the second
line of (\ref{4.8}):
\begin{eqnarray*}
& -8\lambda zv_{z}\cdot 4\lambda ^{2}\left( z^{2}-z/\left( 2\lambda \right)
\right) v=\left( -16\lambda ^{3}\left( z^{3}-z^{2}/\left( 2\lambda \right)
\right) v^{2}\right) _{z}+ \\
&\hspace{3 cm} +48\lambda ^{3}z^{2}\left( 1-\frac{2}{3\lambda z}\right)
v^{2}.
\end{eqnarray*}
Thus, taking
\begin{equation}
\lambda _{0,1}=\frac{9}{3a},  \label{4.10}
\end{equation}%
we obtain
\begin{eqnarray}
& -8\lambda zv_{z}\cdot 4\lambda ^{2}\left( z^{2}-z/\left( 2\lambda \right)
\right) v\geq 43\lambda ^{3}z^{2}u^{2}e^{2\lambda z^{2}}+  \label{4.11} \\
&\hspace{0.1 cm} +\left( -16\lambda ^{3}\left( z^{3}-z^{2}/\left( 2\lambda
\right) \right) u^{2}e^{2\lambda z^{2}}\right) _{z},\forall \lambda \geq
\lambda _{0,1}.  \nonumber
\end{eqnarray}

Summing up (\ref{4.9}) and (\ref{4.11}) and taking into account (\ref{4.8}),
we obtain
\begin{equation}
\hspace{-1.5 cm}
\left.
\begin{array}{c}
\left( \Delta u\right) ^{2}e^{2\lambda z^{2}}\geq -4\lambda
\sum\limits_{i=1}^{n-1}u_{x_{i}}^{2}e^{2\lambda z^{2}}+43\lambda
^{3}z^{2}u^{2}e^{2\lambda z^{2}}+ \\
\left( -4\lambda z\left( u_{z}+2\lambda zu\right) ^{2}e^{2\lambda
z^{2}}+4\lambda z\sum\limits_{i=1}^{n-1}u_{x_{i}}^{2}e^{2\lambda
z^{2}}-16\lambda ^{3}\left( z^{3}- \frac{z^{2}}{2\lambda} \right)
u^{2}e^{2\lambda z^{2}}\right) _{z} \\
+\sum\limits_{i=1}^{n-1}\left( -8\lambda z\left( u_{z}+2\lambda zu\right)
u_{x_{i}}e^{2\lambda z^{2}}\right) _{x_{i}}.%
\end{array}
\right.  \label{4.12}
\end{equation}

\textbf{Step 3. }Consider
\begin{equation*}
\hspace{-2 cm}
\left.
\begin{array}{c}
-\Delta u\cdot ue^{2\lambda
z^{2}}=-\sum\limits_{i=1}^{n-1}u_{x_{i}x_{i}}ue^{2\lambda
z^{2}}-u_{zz}ue^{2\lambda z^{2}}= \\
=\sum\limits_{i=1}^{n-1}\left( -u_{x_{i}}ue^{2\lambda z^{2}}\right)
_{x_{i}}+\sum\limits_{i=1}^{n-1}u_{x_{i}}^{2}e^{2\lambda z^{2}}+\left(
-u_{z}ue^{2\lambda z^{2}}\right) _{z}+u_{z}^{2}e^{2\lambda z^{2}}+4\lambda
zu_{z}ue^{2\lambda z^{2}}= \\
\left( \sum\limits_{i=1}^{n-1}u_{x_{i}}^{2}+u_{z}^{2}\right) e^{2\lambda
z^{2}}+\left( -u_{z}ue^{2\lambda z^{2}}+2\lambda zu^{2}e^{2\lambda
z^{2}}\right) _{z}-8\lambda ^{2}z^{2}\left( 1+\frac{1}{2\lambda z^{2}}
\right) u^{2}e^{2\lambda z^{2}}+ \\
+\sum\limits_{i=1}^{n-1}\left( -u_{x_{i}}ue^{2\lambda z^{2}}\right) _{x_{i}}.%
\end{array}
\right.
\end{equation*}
Thus,
\begin{equation}
\hspace{-1 cm}
\left.
\begin{array}{c}
-\Delta u\cdot ue^{2\lambda z^{2}}=\left(
\sum\limits_{i=1}^{n-1}u_{x_{i}}^{2}+u_{z}^{2}\right) e^{2\lambda
z^{2}}-8\lambda ^{2}z^{2}\left( 1+\frac{1}{2\lambda z^{2}}\right)
u^{2}e^{2\lambda z^{2}}+ \\
+\left( -u_{z}ue^{2\lambda z^{2}}+2\lambda zu^{2}e^{2\lambda z^{2}}\right)
_{z}+\sum\limits_{i=1}^{n-1}\left( -u_{x_{i}}ue^{2\lambda z^{2}}\right)
_{x_{i}}.%
\end{array}
\right.  \label{4.13}
\end{equation}
Multiply (\ref{4.13}) by $5\lambda $ and sum up with (\ref{4.12}). We obtain
\begin{equation*}
\hspace{-2 cm}
\left.
\begin{array}{c}
-5\lambda \Delta u\cdot ue^{2\lambda z^{2}}+\left( \Delta u\right)
^{2}e^{2\lambda z^{2}}\geq \lambda \left( \nabla u\right) ^{2}e^{2\lambda
z^{2}}+ \\
+43\lambda ^{3}z^{2}\left( 1-\frac{40}{43}\left( 1+\frac{1}{2\lambda z^{2}}
\right) u^{2}e^{2\lambda z^{2}}\right) + \\
\left( -4\lambda z\left( u_{z}+2\lambda zu\right) ^{2}e^{2\lambda
z^{2}}+4\lambda z\sum\limits_{i=1}^{n-1}u_{x_{i}}^{2}e^{2\lambda
z^{2}}-16\lambda ^{3}\left( z^{3}- \frac{z^{2}}{2\lambda} \right)
u^{2}e^{2\lambda z^{2}}\right) _{z} \\
+\left( -5u_{z}ue^{2\lambda z^{2}}+10\lambda zu^{2}e^{2\lambda z^{2}}\right)
_{z}+ \\
+\sum\limits_{i=1}^{n-1}\left( -8\lambda z\left( u_{z}+2\lambda zu\right)
u_{x_{i}}e^{2\lambda z^{2}}+\sum\limits_{i=1}^{n-1}\left( -u_{x_{i}}u\right)
e^{2\lambda z^{2}}\right) _{x_{i}},\forall \lambda \geq \lambda _{0,1}.%
\end{array}
\right.
\end{equation*}
Integrating this over the domain $\Omega $ and using Gauss formula, (\ref%
{2.1})-(\ref{2.4}), (\ref{4.05}) and (\ref{4.500}), we obtain that there
exists a sufficiently large number $\lambda _{0}=\lambda _{0}\left( \Omega
\right) \geq \lambda _{01}$ and a number $C=C\left( \Omega \right) >0$ such
that
\begin{equation}
\hspace{-1.5 cm}
\eqalign{
&\hspace{0 cm} \int\limits_{\Omega }\left( -5\lambda \Delta u\cdot
ue^{2\lambda z^{2}}\right) dx+\int\limits_{\Omega }\left( \Delta u\right)
^{2}e^{2\lambda z^{2}}dx C\int\limits_{\Omega }\left( \lambda \left( \nabla
u\right) ^{2}+\lambda ^{3}u^{2}\right) e^{2\lambda z^{2}}dx \\
&-C\lambda ^{3}\left( \left\Vert u\right\Vert _{H^{1}\left( \partial
_{2}\Omega \right) }^{2}+\left\Vert \partial _{z}u\right\Vert _{L_{2}\left(
\partial _{2}\Omega \right) }^{2}\right) e^{2\lambda b^{2}}, \ \forall
\lambda \geq \lambda _{0}, \ \forall u\in C^{2}\left( \overline{\Omega }%
\right) \cap H_{0}^{2}\left( \Omega \right) .
}
\label{4.14}
\end{equation}
By Cauchy-Schwarz inequality $5\lambda \Delta u\cdot ue^{2\lambda z^{2}}\leq
2.5\left( \Delta u\right) ^{2}e^{2\lambda z^{2}}+2.5\lambda
^{2}u^{2}e^{2\lambda z^{2}}.$ Substituting this in (\ref{4.14}) and using
density arguments, we obtain the target estimate (\ref{4.5}). $\ \square $

\subsection{Proof of Theorem 4.1}

\label{sec:4.2}

Denote
\begin{equation}
p\left( x,\alpha \right) =p_{\varepsilon ,1}\left( x,\alpha \right)
-p_{\varepsilon ,2}\left( x,\alpha \right) ,q\left( x,\alpha \right)
=q_{\varepsilon ,1}\left( x,\alpha \right) -q_{\varepsilon ,2}\left(
x,\alpha \right) .  \label{4.15}
\end{equation}%
Then by (\ref{4.2}) and (\ref{4.05})%
\begin{equation}
p,q\in H_{0}^{2}\left( \Omega \right) \times C\left[ -d,d\right] .
\label{4.16}
\end{equation}%
Let the function $f:\mathbb{R}^{k}\rightarrow \mathbb{R}$, $f\in C^{1}\left(
\mathbb{R}^{k}\right) ,k\geq 1.$ It is well known that the following formula
is valid%
\begin{equation}
f\left( y_{1}\right) -f\left( y_{2}\right) =\widetilde{f}\left(
y_{1},y_{2}\right) \circ \left( y_{1}-y_{2}\right) ,  \label{4.17}
\end{equation}%
where \textquotedblleft $\circ "$ is the scalar product in $\mathbb{R}^{k}$
and the vector function $\widetilde{f}\left( y_{1},y_{2}\right) $ is such
that
\begin{equation}
\left\vert \widetilde{f}\left( y_{1},y_{2}\right) \right\vert \leq \max_{
\mathbb{R}^{k}}\left\vert \nabla f\left( y\right) \right\vert .  \label{4.18}
\end{equation}%
Consider the differences $L_{1}\left( p_{\varepsilon ,1},q_{\varepsilon
1}\right) -L_{1}\left( p_{\varepsilon ,2},q_{\varepsilon 2}\right) $ and $%
L_{2}\left( p_{\varepsilon ,1},q_{\varepsilon 1}\right) -L_{2}\left(
p_{\varepsilon ,2},q_{\varepsilon 2}\right) .$ Then, using (\ref{3.7})-(\ref%
{3.9}) and (\ref{4.15})-(\ref{4.18}), we obtain two integral differential
inequalities for $x\in \Omega ,\alpha \in \left( -d,d\right) $
\begin{eqnarray}
&\varepsilon \left\vert \Delta p\right\vert \leq C_{2}\left( \left\vert
\nabla p\right\vert +\left\vert \nabla q\right\vert +\left\vert p\right\vert
+\left\vert q\right\vert +\int\limits_{\Phi _{d}}\left( \left\vert
p\right\vert +\left\vert q\right\vert \right) \left( x,\beta \right) d\beta
\right) ,  \label{4.19} \\
& \varepsilon \left\vert \Delta q\right\vert \leq C_{2}\left( \left\vert
\nabla p\right\vert +\left\vert \nabla q\right\vert +\left\vert p\right\vert
+\left\vert q\right\vert +\int\limits_{\Phi _{d}}\left( \left\vert
p\right\vert +\left\vert q\right\vert \right) \left( x,\beta \right) d\beta
\right) ,  \label{4.20} \\
& \hspace{0.3 cm} p\mid _{\partial _{2}\Omega }=\left( p_{\varepsilon
,1}-p_{\varepsilon ,2}\right) \mid _{\partial _{2}\Omega },\partial
_{z}p\mid _{\partial _{2}\Omega }=\left( \partial _{z}p_{\varepsilon
,1}-\partial _{z}p_{\varepsilon ,2}\right) \mid _{\partial _{2}\Omega },
\label{4.21} \\
& \hspace{0.3 cm} q\mid _{\partial _{2}\Omega }=\left( q_{\varepsilon
,1}-q_{\varepsilon ,2}\right) \mid _{\partial _{2}\Omega },\partial
_{z}q\mid _{\partial _{2}\Omega }=\left( \partial _{z}q_{\varepsilon
,1}-\partial _{z}q_{\varepsilon ,2}\right) \mid _{\partial _{2}\Omega }.
\label{4.22}
\end{eqnarray}%
Here and below in this proof $C_{2}>0$ denotes different positive constant
depending on the same parameters as ones listed in (\ref{4.02}), except of $%
\varepsilon .$

Square both sides of each of equations (\ref{4.19}), (\ref{4.20}), then
multiply by $e^{2\lambda z^{2}}$and then integrate over the domain $\Omega ,$
assuming that $\lambda \geq \lambda _{0},$where $\lambda _{0}$ is the
parameter of Theorem 4.2. And then sum up two resulting inequalities. Using
Cauchy-Schwarz inequality, we obtain
\begin{equation}
\left.
\begin{array}{c}
C_{2}\int\limits_{\Omega }\left[ \left\vert \nabla p\right\vert
^{2}+\left\vert \nabla q\right\vert ^{2}+p^{2}+q^{2}+\int\limits_{\Phi
_{d}}\left( p^{2}+q^{2}\right) \left( x,\beta \right) d\beta \right] \geq \\
\geq \varepsilon ^{2}\int\limits_{\Omega }\left( \left( \Delta p\right)
^{2}+\left( \Delta q\right) ^{2}\right) e^{2\lambda z^{2}}dx,\ \alpha
\in \left( -d,d\right) .%
\end{array}
\right.  \label{4.23}
\end{equation}
Recalling (\ref{4.16}) and applying the Carleman estimate (\ref{4.5}) to the
second line of (\ref{4.23}), we obtain for all $\lambda \geq \lambda _{0}$
\begin{equation}
\hspace{-2 cm}
\left.
\begin{array}{c}
C_{2}\int\limits_{\Omega }\left[ \left\vert \nabla p\right\vert
^{2}+\left\vert \nabla q\right\vert ^{2}+p^{2}+q^{2}+\int\limits_{\Phi
_{d}}\left( p^{2}+q^{2}\right) \left( x,\beta \right) d\beta \right] \geq \\
\geq C\varepsilon ^{2}\int\limits_{\Omega }\left[ \lambda \left( \left(
\nabla p\right) ^{2}+\left( \nabla q\right) ^{2}\right) +\lambda ^{3}\left(
p^{2}+q^{2}\right) \right] e^{2\lambda z^{2}}dx- \\
-C\varepsilon ^{2}\lambda ^{3}\left( \left\Vert p\right\Vert _{H^{1}\left(
\partial _{2}\Omega \right) }^{2}+\left\Vert q\right\Vert _{H^{1}\left(
\partial _{2}\Omega \right) }^{2}+\left\Vert \partial _{z}p\right\Vert
_{L_{2}\left( \partial _{2}\Omega \right) }^{2}+\left\Vert \partial
_{z}q\right\Vert _{L_{2}\left( \partial _{2}\Omega \right) }^{2}\right)
e^{2\lambda b^{2}}.%
\end{array}
\right.  \label{4.24}
\end{equation}
Choose
\begin{equation*}
\lambda _{1}=\lambda _{1}\left( \Omega ,d,\varepsilon ,\left\Vert
G\right\Vert _{C\left( \overline{\Omega }\right) \times C^{1}\left[ -d,d %
\right] ^{2}},\left\Vert \mu _{a}\right\Vert _{C\left( \overline{\Omega }
\right) },\left\Vert \mu _{s}\right\Vert _{C\left( \overline{\Omega }\right)
},M\right) \geq \lambda _{0}\geq 1
\end{equation*}%
such that $C\varepsilon ^{2}\lambda _{1}\geq 2C_{2}$ and then set in (\ref%
{4.24}) $\lambda \geq \lambda _{1}.$ We obtain with the constant $C_{1}>0$
as in (\ref{4.02})
\begin{equation}
\hspace{-2 cm}
\left.
\begin{array}{c}
C_{1}\lambda ^{3}\left( \left\Vert p\right\Vert _{H^{1}\left( \partial
_{2}\Omega \right) }^{2}+\left\Vert q\right\Vert _{H^{1}\left( \partial
_{2}\Omega \right) }^{2}+\left\Vert \partial _{z}p\right\Vert _{L_{2}\left(
\partial _{2}\Omega \right) }^{2}+\left\Vert \partial _{z}q\right\Vert
_{L_{2}\left( \partial _{2}\Omega \right) }^{2}\right) e^{2\lambda b^{2}}+
\\
+C_{1}\int\limits_{\Phi _{d}}\left( p^{2}+q^{2}\right) \left( x,\beta
\right) d\beta \geq \\
\geq \int\limits_{\Omega }\left[ \lambda \left( \left( \nabla p\right)
^{2}+\left( \nabla q\right) ^{2}\right) +\lambda ^{3}\left(
p^{2}+q^{2}\right) \right] e^{2\lambda z^{2}}dx,\ \alpha \in \left(
-d,d\right) .%
\end{array}
\right.  \label{4.25}
\end{equation}
Integrate both sides of (\ref{4.25}) with respect to $\alpha \in \left(
-d,d\right) .$ Then choose $\lambda _{2}\geq \lambda _{1}$ such that $%
\lambda _{2}^{3}\geq 2C_{1}.$ Then set $\lambda =\lambda _{2}.$ We obtain
\begin{equation*}
\left.
\begin{array}{c}
\left\Vert \left( p,q\right) \right\Vert _{H_{2}^{1}\left( \Omega \right)
\times L_{2,2}\left( -d,d\right) }^{2}\leq \\
\leq C_{1}\left( \left\Vert p\right\Vert _{H^{1}\left( \partial _{2}\Omega
\right) \times L_{2}\left( -d,d\right) }^{2}+\left\Vert q\right\Vert
_{H^{1}\left( \partial _{2}\Omega \right) \times L_{2}\left( -d,d\right)
}^{2}\right) + \\
+C_{1}\left( \left\Vert \partial _{z}p\right\Vert _{L_{2}\left( \partial
_{2}\Omega \right) \times L_{2}\left( -d,d\right) }^{2}+\left\Vert \partial
_{z}q\right\Vert _{L_{2}\left( \partial _{2}\Omega \right) \times
L_{2}\left( -d,d\right) }^{2}\right) ,%
\end{array}
\right.
\end{equation*}
which is equivalent with (\ref{4.3}). Estimate (\ref{4.4}) obviously follows
from (\ref{4.3}). \ $\square $

\section{Convergence Analysis for the Convexification for BVP (\protect\ref%
{3.8})-(\protect\ref{3.11})}

\label{sec:5}

Consider an integer $k$,
\begin{equation}
k>\left[ \frac{n+1}{2}\right] +1,  \label{5.1}
\end{equation}%
where $\left[ (n+1)/2\right] $ is the largest integer, which does not exceed
$(n+1)/2.$ Then by embedding theorem and (\ref{5.1})%
\begin{equation}
H^{k}\left( \Omega \right) \subset C^{1}\left( \overline{\Omega }\right)
,\left\Vert f\right\Vert _{C^{1}\left( \overline{\Omega }\right) }\leq
C\left\Vert f\right\Vert _{H^{k}\left( \Omega \right) },\ \forall f\in
H^{k}\left( \Omega \right) .  \label{5.2}
\end{equation}%
Introduce the space $S,$
\begin{equation}
\hspace{-1 cm}
S=\left\{
\begin{array}{c}
\left( p,q\right) \left( x,\alpha \right) : \\
\left\Vert \left( p,q\right) \right\Vert _{S}^{2}=\int\limits_{-d}^{d}\left(
\left\Vert p\left( x,\alpha \right) \right\Vert _{H^{k}\left( \Omega \right)
}^{2}+\left\Vert q\left( x,\alpha \right) \right\Vert _{H^{k}\left( \Omega
\right) }^{2}\right) d\alpha <\infty%
\end{array}
\right\} ,  \label{5.02}
\end{equation}

Let $R>0$ be an arbitrary number. Consider the set of 2D vector functions $%
B\left( g_{1},g_{2},g_{3},g_{4},R\right) $ defined as%
\begin{equation}
B\left( g_{1},g_{2},g_{3},g_{4},R\right) =\left\{
\begin{array}{c}
\left( p,q\right) \in S:p\mid _{\partial \Omega }=g_{1},q\mid _{\partial
\Omega }=g_{2}, \\
p_{z}\mid _{\partial _{2}\Omega }=g_{3},q_{z}\mid _{\partial _{2}\Omega
}=g_{4}, \\
\left\Vert \left( p,q\right) \right\Vert _{S}<R%
\end{array}
\right\} .  \label{5.3}
\end{equation}

To solve BVP (\ref{3.8})-(\ref{3.11}), we consider

\textbf{Minimization Problem.} \emph{Minimize the functional }$J_{\lambda
,\beta }\left( p_{\varepsilon },q_{\varepsilon }\right) $\emph{\ on the set }%
$B( g_{1}$, $g_{2},g_{3},g_{4},R) ,$\emph{\ where}
\begin{equation}
\hspace{-2 cm}
\eqalign{
J_{\lambda ,\gamma }\left( p,q\right) = e^{-2\lambda
b^{2}}\int\limits_{-d}^{d}\int\limits_{\Omega }\left\{ \left[ L_{1}\left(
p,q\right) \right] ^{2}+\left[ L_{2}\left( p,q\right) \right] ^{2}\right\}
e^{2\lambda z^{2}}dxd\alpha + \gamma \left\Vert \left( p,q\right)
\right\Vert _{S}^{2}.
}
\label{5.4}
\end{equation}

In (\ref{5.4}) $\gamma \in \left( 0,1\right) $ is the regularization
parameter and the multiplier $e^{-2\lambda b^{2}}$ is introduced to balance
two terms in the right hand side of (\ref{5.4}). Indeed,
\begin{equation}
\max_{\overline{\Omega }}e^{2\lambda z^{2}}=e^{2\lambda b^{2}},\
\min_{ \overline{\Omega }}e^{2\lambda z^{2}}=e^{2\lambda a^{2}}.
\label{5.40}
\end{equation}%
Let $S$ be the space defined in (\ref{5.02}). We introduce the subspace $%
S_{0}\subset S$ as:%
\begin{equation}
S_{0}=\left\{
\begin{array}{c}
\left( p,q\right) \left( x,\alpha \right) \in S: \\
\left( p,q\right) \mid _{\partial \Omega }=\left( \partial _{z}p,\partial
_{z}q\right) \mid _{\partial _{2}\Omega }=0%
\end{array}
\right\} .  \label{5.5}
\end{equation}

\subsection{Global strict convexity}

\label{sec:5.1}

Below in section 5 $C_{3}>0$ denotes different constants, all of which
depend on the following parameters:%
\begin{equation}
C_{3}=C_{3}\left( \Omega ,d,\varepsilon ,R\right) >0.  \label{5.61}
\end{equation}%
Recall the definition (\ref{5.002})\emph{\ }of the space $H_{2}^{1}\left(
\Omega \right) \times L_{2,2}\left( -d,d\right) .$

\textbf{Theorem 5.1} (global strict convexity). \emph{Denote }%
\begin{equation}
v=\left( p,q\right) .  \label{5.610}
\end{equation}%
\emph{\ For any }$\lambda >0$\emph{\ functional (\ref{5.4}) has the Fr\'{e}
chet derivative }$J_{\lambda ,\gamma }^{\prime }\left( v\right) \in S_{0}$%
\emph{\ at every point }$v\in \overline{B\left(
g_{1},g_{2},g_{3},g_{4},R\right) }.$\emph{\ This derivative is Lipschitz
continuous on }$\overline{B\left( g_{1},g_{2},g_{3},g_{4},R\right) },$ \emph{%
\ \ i.e. there exists a constant }$D>0$\emph{\ such that}
\begin{equation}
\hspace{-1 cm}
\left\Vert J_{\lambda ,\gamma }^{\prime }\left( v_{1}\right) -J_{\lambda
,\gamma }^{\prime }\left( v_{2}\right) \right\Vert _{S}\leq D\left\Vert
v_{1}-v_{2}\right\Vert _{S},\forall v_{1},v_{2}\in \overline{B\left(
g_{1},g_{2},g_{3},g_{4},R\right) }.  \label{5.6100}
\end{equation}%
\emph{Let }$\lambda _{0}\geq 1$\emph{\ be the number of Theorem 4.2. There
exists a sufficiently large number }%
\begin{equation}
\lambda _{3}=\lambda _{3}\left( \Omega ,d,\varepsilon ,R\right) \geq \lambda
_{0}\geq 1  \label{5.7}
\end{equation}%
\emph{such that for every }$\lambda \geq \lambda _{3}$\emph{\ the functional
}$J_{\lambda ,\gamma }\left( v\right) $\emph{\ is strictly convex on the set
}$\overline{B\left( g_{1},g_{2},g_{3},g_{4},R\right) },$\emph{\ i.e. for all
}$v_{1},v_{2}\in \overline{B\left( g_{1},g_{2},g_{3},g_{4},R\right) }$\emph{%
\ \ the following inequality holds:}
\begin{equation}
\left.
\begin{array}{c}
J_{\lambda ,\gamma }\left( v_{2}\right) -J_{\lambda ,\gamma }\left(
v_{1}\right) -J_{\lambda ,\gamma }^{\prime }\left( v_{1}\right) \left(
v_{2}-v_{1}\right) \geq \\
\geq C_{3}e^{-2\lambda \left( b^{2}-a^{2}\right) }\left\Vert
v_{2}-v_{1}\right\Vert _{H_{2}^{1}\left( \Omega \right) \times L_{2,2}\left(
-d,d\right) }^{2}+\gamma \left\Vert v_{2}-v_{1}\right\Vert _{S}^{2}.%
\end{array}
\right.  \label{5.8}
\end{equation}
\emph{Furthermore, for every }$\lambda \geq \lambda _{3}$\emph{\ there
exists unique minimizer }$v_{\min ,\lambda }$\emph{\ of the functional }$%
J_{\lambda ,\gamma }\left( v\right) $\emph{\ on the set }$\overline{B\left(
g_{1},g_{2},g_{3},g_{4},R\right) }$\emph{\ and the following inequality
holds:}%
\begin{equation}
J_{\lambda ,\gamma }^{\prime }\left( v_{\min ,\lambda }\right) \left(
v-v_{\min ,\lambda }\right) \geq 0,\ \forall v\in \overline{B\left(
g_{1},g_{2},g_{3},g_{4},R\right) }.  \label{5.10}
\end{equation}

\textbf{Proof}. Let $v_{1}=\left( p_{1},q_{1}\right) $ and $v_{2}=\left(
p_{2},q_{2}\right) $ be two arbitrary points of the set $\overline{B\left(
g_{1},g_{2},g_{3},g_{4},R\right) }.$ Denote
\begin{equation}
h=\left( h_{1},h_{2}\right) =v_{2}-v_{1}=\left(
p_{2}-p_{1},q_{2}-q_{1}\right) .  \label{5.11}
\end{equation}%
By (\ref{5.2}), (\ref{5.02}), (\ref{5.3}), (\ref{5.5}), (\ref{5.61}) and (%
\ref{5.11})
\begin{eqnarray}
&\hspace{3 cm} h\in S_{0},  \label{5.13} \\
& \left\Vert h\right\Vert _{S}\leq 2R,\left\Vert h_{i}\right\Vert
_{C^{1}\left( \overline{\Omega }\right) \times C\left[ -d,d\right] }\leq
C_{3},i=1,2.  \label{5.130}
\end{eqnarray}%
By (\ref{3.8}) and (\ref{5.11})
\begin{equation}
\hspace{-2 cm}
\left.
\begin{array}{c}
L_{1}\left( p_{2},q_{2}\right) =L_{1}\left( p_{1}+h_{1},q_{1}+h_{2}\right) =
\\
=L_{1}\left( p_{1},q_{1}\right) +\left[ -\varepsilon \Delta h_{1}+\nu (x
\mathbf{,}\alpha )\cdot \nabla _{x}h_{1}(x,\alpha )+\partial _{\alpha }\nu
(x \mathbf{,}\alpha )\cdot \nabla _{x}h_{1}(x,\alpha )\right] + \\
+h_{2}e^{-\left( p_{1}+h_{1}\right) \left( x,\alpha \right) }\mu
_{s}(x)\int\limits_{\Phi _{d}}G(x,\alpha ,\beta )e^{\left(
p_{1}+h_{1}\right) \left( x,\beta \right) }d\beta - \\
-e^{-\left( p_{1}+h_{1}\right) \left( x,\alpha \right) }\mu
_{s}(x)\int\limits_{\Phi _{d}}\partial _{\alpha }G(x,\alpha ,\beta
)e^{\left( p_{1}+h_{1}\right) \left( x,\beta \right) }d\beta + \\
+e^{-p_{1}\left( x,\alpha \right) }\mu _{s}(x)\int\limits_{\Phi
_{d}}\partial _{\alpha }G(x,\alpha ,\beta )e^{p_{1}\left( x,\beta \right)
}d\beta .%
\end{array}
\right.  \label{5.14}
\end{equation}
Consider now the nonlinear term in the third line of (\ref{5.14}). By Taylor
formula%
\begin{eqnarray}
& e^{-\left( p_{1}+h_{1}\right) \left( x,\alpha \right) }=e^{-p_{1}\left(
x,\alpha \right) }-e^{-p_{1}\left( x,\alpha \right) }h_{1}\left( x,\alpha
\right) +\frac{1}{2}e^{-\xi }h_{1}^{2}\left( x,\alpha \right) ,
\label{5.140} \\
&\hspace{0.5 cm} e^{\left( p_{1}+h_{1}\right) \left( x,\beta \right)
}=e^{p_{1}\left( x,\beta \right) }+e^{p_{1}\left( x,\beta \right)
}h_{1}\left( x,\beta \right) +\frac{ 1}{2}e^{\eta }h_{1}^{2}\left( x,\beta
\right) ,  \label{5.141}
\end{eqnarray}%
where $\xi $ is a point between $p_{1}\left( x,\alpha \right) $ and $%
p_{1}\left( x,\alpha \right) +h_{1}\left( x,\alpha \right) $ and $\eta $ is
a point between $p_{1}\left( x,\beta \right) $ and $p_{1}\left( x,\beta
\right) +h_{1}\left( x,\beta \right) .$ Hence, using (\ref{5.130}), (\ref%
{5.140}) and (\ref{5.141}), we obtain%
\begin{equation}
\left.
\begin{array}{c}
h_{2}e^{-\left( p_{1}+h_{1}\right) \left( x,\alpha \right) }\mu
_{s}(x)\int\limits_{\Phi _{d}}G(x,\alpha ,\beta )e^{\left(
p_{1}+h_{1}\right) \left( x,\beta \right) }d\beta = \\
=h_{2}e^{-p_{1}\left( x,\alpha \right) }\mu _{s}(x)\int\limits_{\Phi
_{d}}G(x,\alpha ,\beta )e^{p_{1}\left( x,\beta \right) }d\beta +Q_{1}\left(
h,x,\alpha \right) , \\
\left\vert Q_{1}\left( h,x,\alpha \right) \right\vert \leq C_{3}\left(
h^{2}+\int\limits_{\Phi _{d}}h^{2}\left( x,\beta \right) d\beta \right) , \\
x\in \Omega ,\alpha \in \left( -d,d\right) .%
\end{array}
\right.  \label{5.15}
\end{equation}%
And also%
\begin{equation}
\left.
\begin{array}{c}
-e^{-\left( p_{1}+h_{1}\right) \left( x,\alpha \right) }\mu
_{s}(x)\int\limits_{\Phi _{d}}\partial _{\alpha }G(x,\alpha ,\beta
)e^{\left( p_{1}+h_{1}\right) \left( x,\beta \right) }d\beta + \\
+e^{-p_{1}\left( x,\alpha \right) }\mu _{s}(x)\int\limits_{\Phi
_{d}}\partial _{\alpha }G(x,\alpha ,\beta )e^{p_{1}\left( x,\beta \right)
}d\beta = \\
=h_{1}\left( x,\alpha \right) e^{-p_{1}\left( x,\alpha \right) }\mu
_{s}(x)\int\limits_{\Phi _{d}}\partial _{\alpha }G(x,\alpha ,\beta
)e^{p_{1}\left( x,\beta \right) }d\beta - \\
-e^{-p_{1}\left( x,\alpha \right) }\mu _{s}(x)\int\limits_{\Phi
_{d}}\partial _{\alpha }G(x,\alpha ,\beta )e^{p_{1}\left( x,\beta \right)
}h_{1}\left( x,\beta \right) d\beta +Q_{2}\left( h,x,\alpha \right) , \\
\left\vert Q_{2}\left( h,x,\alpha \right) \right\vert \leq C_{3}\left(
h^{2}+\int\limits_{\Phi _{d}}h^{2}\left( x,\beta \right) d\beta \right) , \\
x\in \Omega ,\alpha \in \left( -d,d\right) .%
\end{array}
\right.  \label{5.16}
\end{equation}

Thus, by (\ref{5.14}), (\ref{5.15}) and (\ref{5.16})%
\begin{equation}
L_{1}\left( p_{2},q_{2}\right) =L_{1}\left( p_{1},q_{1}\right) +L_{1,\mbox{
lin}}\left( p_{1},q_{1,}h\right) +L_{1,\mbox{nonlin}}\left(
p_{1},q_{1,}h\right) ,  \label{5.17}
\end{equation}%
where $L_{1,\mbox{lin}}\left( p_{1},q_{1,}h\right) $ depends linearly on $h$
and $L_{1,\mbox{nonlin}}\left( p_{1},q_{1,}h\right) $ depends nonlinearly on
$h.$ Also, for $x\in \Omega ,\alpha \in \left( -d,d\right) $%
\begin{eqnarray}
&\hspace{-0.5 cm} L_{1,\mbox{lin}}\left( p_{1},q_{1,}h\right) \left( x,\alpha
\right) =-\varepsilon \Delta h_{1}\left( x,\alpha \right) +\widehat{L}_{1,%
\mbox{lin} }\left( p_{1},q_{1,}h\right) \left( x,\alpha \right) ,
\label{5.18} \\
&\hspace{-1.5 cm}\left\vert \widehat{L}_{1,\mbox{lin}}\left( p_{1},q_{1,}h\right) \left(
x,\alpha \right) \right\vert \leq C_{3}\left( \left\vert \nabla
_{x}h_{1}\right\vert +\left\vert \nabla _{x}h_{2}\right\vert +\left\vert
h\right\vert +\int\limits_{\Phi _{d}}\left\vert h\left( x,\beta \right)
\right\vert d\beta \right) ,  \label{5.19} \\
&\hspace{0 cm} \left\vert L_{1,\mbox{nonlin}}\left( p_{1},q_{1,}h\right)
\right\vert \leq C_{3}\left( \left\vert h\left( x,\alpha \right) \right\vert
^{2}+\int\limits_{\Phi _{d}}\left\vert h\left( x,\beta \right) \right\vert
^{2}d\beta \right) .  \label{5.20}
\end{eqnarray}%
Similarly,%
\begin{eqnarray}
& \hspace{0 cm} L_{2}\left( p_{2},q_{2}\right) =L_{2}\left(
p_{1},q_{1}\right) +L_{2,\mbox{ lin}}\left( p_{1},q_{1,}h\right) +L_{2,\mbox{%
nonlin}}\left( p_{1},q_{1,}h\right) ,  \label{5.21} \\
& \hspace{0 cm} L_{2,\mbox{lin}}\left( p_{1},q_{1,}h\right) \left(
x,\alpha \right) =-\varepsilon \Delta h_{2}\left( x,\alpha \right) +\widehat{%
L}_{2,\mbox{lin} }\left( p_{1},q_{1,}h\right) \left( x,\alpha \right) ,
\label{5.22} \\
&\hspace{-1.2 cm} \left\vert \widehat{L}_{2,\mbox{lin}}\left( p_{1},q_{1,}h\right) \left(
x,\alpha \right) \right\vert \leq C_{3}\left( \left\vert \nabla
_{x}h_{1}\right\vert +\left\vert \nabla _{x}h_{2}\right\vert +\left\vert
h\right\vert +\int\limits_{\Phi _{d}}\left\vert h\left( x,\beta \right)
\right\vert d\beta \right) ,  \label{5.23} \\
& \hspace{0 cm} \left\vert L_{2,\mbox{nonlin}}\left( p_{1},q_{1,}h\right)
\right\vert \leq C_{3}\left( \left\vert h\left( x,\alpha \right) \right\vert
^{2}+\int\limits_{\Phi _{d}}\left\vert h\left( x,\beta \right) \right\vert
^{2}d\beta \right) .  \label{5.24}
\end{eqnarray}

Using (\ref{5.4}) and (\ref{5.17})-(\ref{5.24}) and recalling (\ref{5.610}),
we obtain
\begin{equation}
\hspace{-2 cm}
\left.
\begin{array}{c}
J_{\lambda ,\beta }\left( v_{2}\right) -J_{\lambda ,\beta }\left(
v_{1}\right) - \\
-2e^{-2\lambda b^{2}}\int\limits_{-d}^{d}\int\limits_{\Omega }L_{1}\left(
p_{1},q_{1}\right) L_{1,\mbox{lin}}\left( p_{1},q_{1},h\right) e^{2\lambda
z^{2}}dxd\alpha - \\
-2e^{-2\lambda b^{2}}\int\limits_{-d}^{d}\int\limits_{\Omega }L_{2}\left(
p_{1},q_{1}\right) L_{2,\mbox{lin}}\left( p_{1},q_{1},h\right) e^{2\lambda
z^{2}}dxd\alpha -2\gamma \left[ v_{1},h\right] = \\
=e^{-2\lambda b^{2}}\int\limits_{-d}^{d}\int\limits_{\Omega }\left[ \left(
L_{1,\mbox{lin}}\left( p_{1},q_{1},h\right) \right) ^{2}+\left( L_{2,\mbox{
lin}}\left( p_{1},q_{1},h\right) \right) ^{2}\right] e^{2\lambda
z^{2}}dxd\alpha + \\
+e^{-2\lambda b^{2}}\int\limits_{-d}^{d}\int\limits_{\Omega }\left[ \left(
L_{1,\mbox{nonlin}}\left( p_{1},q_{1},h\right) \right) ^{2}+\left( L_{2,
\mbox{nonlin}}\left( p_{1},q_{1},h\right) \right) ^{2}\right] e^{2\lambda
z^{2}}dxd\alpha + \\
+2e^{-2\lambda b^{2}}\int\limits_{-d}^{d}\int\limits_{\Omega }L_{1,\mbox{
nonlin}}\left( p_{1},q_{1},h\right) \left[ L_{1}\left( p_{1},q_{1}\right)
+L_{1,\mbox{lin}}\left( p_{1},q_{1},h\right) \right] e^{2\lambda
z^{2}}dxd\alpha + \\
+2e^{-2\lambda b^{2}}\int\limits_{-d}^{d}\int\limits_{\Omega }L_{2,\mbox{
nonlin}}\left( p_{1},q_{1},h\right) \left[ L_{2}\left( p_{1},q_{1}\right)
+L_{2,\mbox{lin}}\left( p_{1},q_{1},h\right) \right] e^{2\lambda
z^{2}}dxd\alpha + \\
+\gamma \left\Vert h\right\Vert _{S}^{2},%
\end{array}
\right.  \label{5.25}
\end{equation}
where $\left[ ,\right] $ is the scalar product in $S$.

Consider the sum of second and third lines of (\ref{5.25}),%
\begin{equation}
\hspace{-1 cm}
\begin{array}{c}
I\left( p_{1},q_{1},h\right) =-2e^{-2\lambda
b^{2}}\int\limits_{-d}^{d}\int\limits_{\Omega }L_{1}\left(
p_{1},q_{1}\right) L_{1,\mbox{lin}}\left( p_{1},q_{1},h\right) e^{2\lambda
z^{2}}dxd\alpha - \\
-2e^{-2\lambda b^{2}}\int\limits_{-d}^{d}\int\limits_{\Omega }L_{2}\left(
p_{1},q_{1}\right) L_{2,\mbox{lin}}\left( p_{1},q_{1},h\right) e^{2\lambda
z^{2}}dxd\alpha -2\gamma \left[ v_{1},h\right] .%
\end{array}
\label{5.26}
\end{equation}%
It follows from (\ref{5.5}) and (\ref{5.13}) that we can consider $I\left(
p_{1},q_{1},h\right) :S_{0}\rightarrow \mathbb{R}$ as a linear functional of
$h\in S_{0}.$ It follows from (\ref{5.02}), (\ref{5.18}), (\ref{5.19}), (\ref%
{5.22}) and (\ref{5.23}) that $I\left( p_{1},q_{1},h\right) $ is a bounded
functional. Therefore, by Riesz theorem there exists unique point $\widehat{I%
}\left( p_{1},q_{1}\right) \in S_{0}$ such that
\begin{equation}
I\left( p_{1},q_{1},h\right) =\left[ \widehat{I}\left( p_{1},q_{1}\right) ,h %
\right] ,\forall h\in S_{0}.  \label{5.27}
\end{equation}%
It follows from (\ref{5.610}), (\ref{5.11}) and (\ref{5.18})-(\ref{5.27})
that
\begin{equation*}
\lim_{\left\Vert h\right\Vert _{S}\rightarrow 0}\frac{J_{\lambda ,\gamma
}\left( p_{1}+h_{1},q_{1}+h_{2}\right) -J_{\lambda ,\gamma }\left(
p_{1},q_{1}\right) -\left[ \widehat{I}\left( p_{1},q_{1}\right) ,h\right] }{
\left\Vert h\right\Vert _{S}}=0.
\end{equation*}%
Therefore $\widehat{I}\left( p_{1},q_{1}\right) =J_{\lambda ,\gamma
}^{\prime }\left( p_{1},q_{1}\right) $ is the Fr\'{e}chet derivative of the
functional $J_{\lambda ,\beta }\left( p,q\right) $ at the point $\left(
p_{1},q_{1}\right) .$ The Lipschitz continuity property (\ref{5.6100}) is
proven similarly with the proof of Theorem 3.1 of \cite{Bak}. Therefore, we
omit this proof here.

Thus, (\ref{5.25}) can be rewritten as
\begin{equation}
\hspace{-2 cm}
\left.
\begin{array}{c}
J_{\lambda ,\gamma }\left( p_{1}+h_{1},q_{1}+h_{2}\right) -J_{\lambda
,\gamma }\left( p_{1},q_{1}\right) -J_{\lambda ,\gamma }^{\prime }\left(
p_{1},q_{1}\right) \left( h\right) = \\
=e^{-2\lambda b^{2}}\int\limits_{-d}^{d}\int\limits_{\Omega }\left[ \left(
L_{1,\mbox{lin}}\left( p_{1},q_{1},h\right) \right) ^{2}+\left( L_{2,\mbox{
lin}}\left( p_{1},q_{1},h\right) \right) ^{2}\right] e^{2\lambda
z^{2}}dxd\alpha + \\
+e^{-2\lambda b^{2}}\int\limits_{-d}^{d}\int\limits_{\Omega }\left[ \left(
L_{1,\mbox{nonlin}}\left( p_{1},q_{1},h\right) \right) ^{2}+\left( L_{2,
\mbox{nonlin}}\left( p_{1},q_{1},h\right) \right) ^{2}\right] e^{2\lambda
z^{2}}dxd\alpha + \\
+2e^{-2\lambda b^{2}}\int\limits_{-d}^{d}\int\limits_{\Omega }L_{1,\mbox{
nonlin}}\left( p_{1},q_{1},h\right) \left[ L_{1}\left( p_{1},q_{1}\right)
+L_{1,\mbox{lin}}\left( p_{1},q_{1},h\right) \right] e^{2\lambda
z^{2}}dxd\alpha + \\
+2e^{-2\lambda b^{2}}\int\limits_{-d}^{d}\int\limits_{\Omega }L_{2,\mbox{
nonlin}}\left( p_{1},q_{1},h\right) \left[ L_{2}\left( p_{1},q_{1}\right)
+L_{2,\mbox{lin}}\left( p_{1},q_{1},h\right) \right] e^{2\lambda
z^{2}}dxd\alpha + \\
+\gamma \left\Vert h\right\Vert _{S}^{2}.%
\end{array}
\right.  \label{5.28}
\end{equation}

Let $RHS$ denotes the right hand side of (\ref{5.28}). Then using (\ref{5.2}%
), (\ref{5.3}), (\ref{5.130}), (\ref{5.18})-(\ref{5.24}) and Cauchy-Schwarz
inequality, we obtain the following estimate from the below:
\begin{equation}
\left.
\begin{array}{c}
\hspace{1 cm} RHS\geq \frac{\varepsilon ^{2}}{2}e^{-2\lambda
b^{2}}\int\limits_{-d}^{d}\int\limits_{\Omega }\left[ \left( \Delta
h_{1}\right) ^{2}+\left( \Delta h_{2}\right) ^{2}\right] e^{2\lambda
z^{2}}dxd\alpha - \\
-C_{3}e^{-2\lambda b^{2}}\int\limits_{-d}^{d}\int\limits_{\Omega }\left[
\left( \nabla h_{1}\right) ^{2}+\left( \nabla h_{2}\right)
^{2}+h_{1}^{2}+h_{2}^{2}\right] e^{2\lambda z^{2}}dxd\alpha +\gamma
\left\Vert h\right\Vert _{S}^{2}.%
\end{array}
\right.  \label{5.29}
\end{equation}

It follows from (\ref{4.05}), (\ref{5.5}) and (\ref{5.13}) that we can apply
Carleman estimate (\ref{4.5}) of Theorem 4.2 to the right hand side of (\ref%
{5.29}), and the second line of (\ref{4.5}) should be zero in this case. Let
$\lambda _{0}=\lambda _{0}\left( \Omega \right) \geq 1$ be the number, which
was found in Theorem 4.2. We obtain for all $\lambda \geq \lambda _{0}$%
\begin{equation*}
\hspace{-1 cm}
\left.
\begin{array}{c}
RHS\geq C\varepsilon ^{2}e^{-2\lambda
b^{2}}\int\limits_{-d}^{d}\int\limits_{\Omega }\left[ \lambda \left( \left(
\nabla h_{1}\right) ^{2}+\left( \nabla h_{2}\right) ^{2}\right) +\lambda
^{3}\left( h_{1}^{2}+h_{2}^{2}\right) \right] e^{2\lambda z^{2}}dxd\alpha -
\\
\hspace{0.5 cm} -C_{3}e^{-2\lambda
b^{2}}\int\limits_{-d}^{d}\int\limits_{\Omega }\left[ \left( \nabla
h_{1}\right) ^{2}+\left( \nabla h_{2}\right) ^{2}+h_{1}^{2}+h_{2}^{2}\right]
e^{2\lambda z^{2}}dxd\alpha +\gamma \left\Vert h\right\Vert _{S}^{2}.%
\end{array}
\right.
\end{equation*}
Choose the number $\lambda _{3}\geq \lambda _{0}$ depending on the same
parameters as the ones listed in (\ref{5.7}) and such that $C\varepsilon
^{2}\lambda _{3}\geq 2C_{3}.$ Then we obtain
\begin{equation*}
\hspace{-2 cm}
RHS\geq C_{3}e^{-2\lambda b^{2}}\int\limits_{-d}^{d}\int\limits_{\Omega }
\left[ \left( \nabla h_{1}\right) ^{2}+\left( \nabla h_{2}\right)
^{2}+h_{1}^{2}+h_{2}^{2}\right] e^{2\lambda z^{2}}dxd\alpha +\gamma
\left\Vert h\right\Vert _{S}^{2},\ \forall \lambda \geq \lambda _{2}.
\end{equation*}%
This, (\ref{5.40}) and (\ref{5.28}) imply
\begin{eqnarray*}
& J_{\lambda ,\gamma }\left( p_{1}+h_{1},q_{1}+h_{2}\right) -J_{\lambda
,\gamma }\left( p_{1},q_{1}\right) -J_{\lambda ,\gamma }^{\prime }\left(
p_{1},q_{1}\right) \left( h\right) \geq \\
&\hspace{0.5 cm} \geq C_{3}e^{-2\lambda \left( b^{2}-a^{2}\right)
}\left\Vert h\right\Vert _{H_{2}^{1}\left( \Omega \right) \times
L_{2,2}\left( -d,d\right) }^{2}+\gamma \left\Vert h\right\Vert _{S}^{2},
\end{eqnarray*}%
which is equivalent with (\ref{5.8}).

Given $\lambda \geq \lambda _{3},$ the existence and uniqueness of the
minimizer $v_{\min ,\lambda }$ of the functional $J_{\lambda ,\gamma }\left(
v\right) $ on the set $\overline{B\left( g_{1},g_{2},g_{3},g_{4},R\right) }$
as well as inequality (\ref{5.10}) follow immediately from a combination of
either Lemma 2.1 with Theorem 2.1 of \cite{Bak} or, equivalently, Lemma
5.2.1 and Theorem 5.2.1 of \cite{KL}. $\ \ \square $

\subsection{The accuracy of the minimizer}

\label{sec:5.2}

In this section we estimate the distance between the minimizer $v_{\min
,\lambda }$ which was found in Theorem 5.1, and the exact solution $v^{\ast
} $ with the noiseless data of BVP (\ref{3.8})-(\ref{3.11}). In accordance
with the regularization theory \cite{T}, we assume that there exists a
solution $v_{\varepsilon }^{\ast }\left( x,\alpha \right) =\left(
p_{\varepsilon }^{\ast },q_{\varepsilon }^{\ast }\right) \left( x,\alpha
\right) $ of BVP (\ref{3.8})-(\ref{3.11}) with the noiseless boundary data $%
g_{1}^{\ast },g_{2}^{\ast },g_{3}^{\ast },g_{4}^{\ast }$ in (\ref{3.10}), (%
\ref{3.11}). By Theorem 4.1 this solution is unique. However, in
applications the data (\ref{3.10}), (\ref{3.11}) are always given with a
noise. Let a small number $\delta \in \left( 0,1\right) $ be the level of
the noise in the data (\ref{3.10}), (\ref{3.11}). We assume that there exist
two vector functions $F\left( x,\alpha \right) =\left( F_{1},F_{2}\right)
\left( x,\alpha \right) \in S$ and $F^{\ast }\left( x,\alpha \right) =\left(
F_{1}^{\ast },F_{2}^{\ast }\right) \left( x,\alpha \right) \in S$ such that%
\begin{equation}
\left.
\begin{array}{c}
F_{1}\mid _{\partial \Omega }=g_{1}\left( x,\alpha \right) ,\partial
_{z}F_{1}\mid _{\partial _{2}\Omega }=g_{3}\left( x,\alpha \right) , \\
F_{1}^{\ast }\mid _{\partial \Omega }=g_{1}^{\ast }\left( x,\alpha \right)
,\partial _{z}F_{1}^{\ast }\mid _{\partial _{2}\Omega }=g_{3}^{\ast }\left(
x,\alpha \right) , \\
F_{2}\mid _{\partial \Omega }=g_{2}\left( x,\alpha \right) ,\partial
_{z}F_{2}\mid _{\partial _{2}\Omega }=g_{4}\left( x,\alpha \right) , \\
F_{2}^{\ast }\mid _{\partial \Omega }=g_{2}^{\ast }\left( x,\alpha \right)
,\partial _{z}F_{2}^{\ast }\mid _{\partial _{2}\Omega }=g_{4}^{\ast }\left(
x,\alpha \right) .%
\end{array}
\right.  \label{5.30}
\end{equation}%
And we also assume that%
\begin{eqnarray}
&\hspace{0.1 cm} \left\Vert v_{\varepsilon }^{\ast }\right\Vert
_{S}<R-C_{3}\delta ,  \label{5.300} \\
& \left\Vert F\right\Vert _{S},\left\Vert F^{\ast }\right\Vert _{S}<R,
\label{5.31} \\
&\hspace{0.3 cm} \left\Vert F-F^{\ast }\right\Vert _{S}<\delta .
\label{5.32}
\end{eqnarray}%
For every vector function $v=\left( p,q\right) \in B\left(
g_{1},g_{2},g_{3},g_{4},R\right) $ consider the difference
\begin{equation}
\widetilde{v}=\left( \widetilde{p},\widetilde{q}\right) =\left(
p-F_{1},q-F_{2}\right) =v-F.  \label{5.33}
\end{equation}%
Also, denote
\begin{equation}
\widetilde{v}_{\varepsilon }^{\ast }=\left( \widetilde{p}_{\varepsilon
}^{\ast },\widetilde{q}_{\varepsilon }^{\ast }\right) =\left( p_{\varepsilon
}^{\ast }-F_{1}^{\ast },q_{\varepsilon }^{\ast }-F_{2}^{\ast }\right)
=v_{\varepsilon }^{\ast }-F^{\ast }.  \label{5.34}
\end{equation}%
Using (\ref{5.3}), (\ref{5.5}) and (\ref{5.300})-(\ref{5.34}), we obtain%
\begin{equation}
\widetilde{v},\widetilde{v}_{\varepsilon }^{\ast }\in B_{0}\left( 2R\right)
=\left\{ v=\left( p,q\right) \in S_{0}:\left\Vert v\right\Vert
_{S}<2R\right\} .  \label{5.35}
\end{equation}

\textbf{Theorem 5.2.} \emph{Let conditions (\ref{5.300})-(\ref{5.32}) and
notations (\ref{5.33}),(\ref{5.34}) hold. Let the regularization parameter }$%
\gamma $ \emph{be} \emph{\ }%
\begin{equation}
\gamma =\gamma \left( \delta \right) =\delta ^{2}.  \label{5.350}
\end{equation}%
\emph{Let }$\lambda _{3}=\lambda _{3}\left( \Omega ,d,\varepsilon ,R\right)
\geq \lambda _{0}\geq 1$\emph{\ be the number in (\ref{5.7}). Consider the
number }$\lambda _{4},$%
\begin{equation}
\lambda _{4}=\lambda _{3}\left( \Omega ,d,\varepsilon ,2R\right) \geq
\lambda _{3}.  \label{5.36}
\end{equation}%
\emph{Let }$\lambda =\lambda _{4}$\emph{\ and let }$v_{\min ,\lambda
_{4}}=\left( p_{\min ,\lambda _{4}},q_{\min ,\lambda _{4}}\right) \in
\overline{B\left( g_{1},g_{2},g_{3},g_{4},R\right) }$\emph{\ be the
minimizer of the functional }$J_{\lambda _{4},\gamma }\left( v\right) $\emph{%
\ \ \ on the set }$\overline{B\left( g_{1},g_{2},g_{3},g_{4},R\right) }.$%
\emph{\ \ Let }$a_{\min ,\lambda _{4}}\left( x\right) $\emph{\ and }$%
a_{\varepsilon }^{\ast }\left( x\right) $\emph{\ be functions }$a\left(
x\right) ,$\emph{\ which are constructed from }$v_{\min ,\lambda
_{4}}=\left( p_{\min ,\lambda _{4}},q_{\min ,\lambda _{4}}\right) $\emph{\
and }$v_{\varepsilon }^{\ast }=\left( p_{\varepsilon }^{\ast
},q_{\varepsilon }^{\ast }\right) $\emph{\ respectively via the right hand
side of formula (\ref{5.7}), in which }$p_{\varepsilon ,\mbox{comp}}$\emph{\
is replaced with }$p_{\min ,\lambda }$\emph{\ and }$p_{\varepsilon }^{\ast }$%
\emph{\ respectively. Then the following accuracy estimates hold:}%
\begin{eqnarray}
&\left\Vert v_{\min ,\lambda _{4}}-v_{\varepsilon }^{\ast }\right\Vert
_{H_{2}^{1}\left( \Omega \right) \times L_{2,2}\left( -d,d\right) }\leq
C_{3}\delta ,  \label{5.37} \\
&\hspace{0.8 cm} \left\Vert a_{\min ,\lambda _{4}}-a_{\varepsilon }^{\ast
}\right\Vert _{L_{2}\left( \Omega \right) }\leq C_{3}\delta .  \label{5.38}
\end{eqnarray}

\textbf{Proof}. Let $B_{0}\left( 2R\right) $ be the set defined in (\ref%
{5.35}). Consider a new functional $I_{\lambda _{4},\gamma }:B_{0}\left(
2R\right) \rightarrow \mathbb{R}$ defined as%
\begin{equation}
I_{\lambda _{4},\gamma }\left( \widetilde{v}\right) =J_{\lambda _{4},\gamma
}\left( \widetilde{v}+F\right) .  \label{5.42}
\end{equation}%
Then Theorem 5.2 is applicable to this functional. Let $V_{\min ,\lambda
_{3}}$ be the minimizer of $I_{\lambda _{4},\gamma }\left( \widetilde{v}%
\right) $ on the set $\overline{B_{0}\left( 2R\right) },$%
\begin{equation}
\min_{\overline{B_{0}\left( 2R\right) }}J_{\lambda _{4},\gamma }\left(
\widetilde{v}+F\right) =J_{\lambda _{4},\gamma }\left( V_{\min ,\lambda
_{3}}+F\right) .  \label{5.420}
\end{equation}%
Since by (\ref{5.35}) both vector functions $V_{\min ,\lambda _{4}},%
\widetilde{v}_{\varepsilon }^{\ast }\in S_{0}\left( 2R\right) ,$ then by an
obvious analog of (\ref{5.8})
\begin{equation}
\left.
\begin{array}{c}
I_{\lambda _{4},\gamma }\left( \widetilde{v}_{\varepsilon }^{\ast }\right)
-I_{\lambda _{4},\gamma }\left( V_{\min ,\lambda }\right) -I_{\lambda
_{4},\gamma }^{\prime }\left( \widetilde{v}_{\min ,\lambda _{4}}\right)
\left( \widetilde{v}_{\varepsilon }^{\ast }-V_{\min ,\lambda _{4}}\right)
\geq \\
\hspace{0.5 cm} \geq C_{3}e^{-2\lambda _{4}\left( b^{2}-a^{2}\right)
}\left\Vert \widetilde{ v }_{\varepsilon }^{\ast }-V_{\min ,\lambda
_{4}}\right\Vert _{H_{2}^{1}\left( \Omega \right) \times L_{2,2}\left(
-d,d\right) }^{2}.%
\end{array}
\right.  \label{5.43}
\end{equation}
By (\ref{5.10})%
\begin{equation*}
-I_{\lambda _{4},\gamma }\left( \widetilde{v}_{\min ,\lambda _{4}}\right)
-I_{\lambda _{4},\gamma }^{\prime }\left( \widetilde{v}_{\min ,\lambda
_{4}}\right) \left( \widetilde{v}_{\varepsilon }^{\ast }-\widetilde{v}_{\min
,\lambda _{4}}\right) \leq 0.
\end{equation*}%
Hence, (\ref{5.43}) implies%
\begin{equation*}
C_{3}e^{2\lambda _{4}\left( b^{2}-a^{2}\right) }I_{\lambda _{4},\gamma
}\left( \widetilde{v}_{\varepsilon }^{\ast }\right) \geq \left\Vert
\widetilde{v}_{\varepsilon }^{\ast }-V_{\min ,\lambda _{4}}\right\Vert
_{H_{2}^{1}\left( \Omega \right) \times L_{2,2}\left( -d,d\right) }^{2}.
\end{equation*}%
Taking into account dependencies (\ref{5.61}) and (\ref{5.36}), we obtain%
\begin{equation}
\left\Vert \widetilde{v}_{\varepsilon }^{\ast }-V_{\min ,\lambda
_{4}}\right\Vert _{H_{2}^{1}\left( \Omega \right) \times L_{2,2}\left(
-d,d\right) }^{2}\leq C_{3}\delta ^{2}.  \label{5.44}
\end{equation}%
By (\ref{5.34}) and (\ref{5.42})
\begin{equation}
I_{\lambda _{4},\gamma }\left( \widetilde{v}_{\varepsilon }^{\ast }\right)
=J_{\lambda _{4},\gamma }\left( \widetilde{v}_{\varepsilon }^{\ast
}+F\right) =J_{\lambda _{4},\gamma }\left( v_{\varepsilon }^{\ast }+\left(
F-F^{\ast }\right) \right) .  \label{5.45}
\end{equation}%
By (\ref{5.4})
\begin{equation}
\hspace{-2 cm}
\eqalign{
J_{\lambda _{4},\gamma }\left( v_{\varepsilon }^{\ast }\right) =
e^{-2\lambda _{4}b^{2}}\int\limits_{-d}^{d}\int\limits_{\Omega }\left\{
\left[ L_{1}\left( v^{\ast }\right) \right] ^{2}+\left[ L_{2}\left( v^{\ast
}\right) \right] ^{2}\right\} e^{2\lambda _{4}z^{2}}dxd\alpha + \gamma
\left\Vert v^{\ast }\right\Vert _{S}^{2}.
}
\label{5.46}
\end{equation}%
Since $v_{\varepsilon }^{\ast }$ is the exact solution of BVP (\ref{3.8})-(%
\ref{3.11}), then $L_{1}\left( v^{\ast }\right) =L_{2}\left( v^{\ast
}\right) =0.$ Hence, (\ref{5.350}) and (\ref{5.46}) imply
\begin{equation}
J_{\lambda _{4},\gamma }\left( v_{\varepsilon }^{\ast }\right) =\gamma
\left\Vert v^{\ast }\right\Vert _{S}^{2}\leq C_{3}\delta ^{2}.  \label{5.47}
\end{equation}%
Next, by (\ref{5.32}), (\ref{5.45})-(\ref{5.47}) and Cauchy-Schwarz
inequality%
\begin{equation*}
I_{\lambda _{4},\gamma }\left( \widetilde{v}_{\varepsilon }^{\ast }\right)
=J_{\lambda _{4},\gamma }\left( v_{\varepsilon }^{\ast }+\left( F-F^{\ast
}\right) \right) \leq C_{3}\left( J_{\lambda _{4},\gamma }\left(
v_{\varepsilon }^{\ast }\right) +\left\Vert F-F^{\ast }\right\Vert
_{S}^{2}\right) \leq C_{3}\delta ^{2}.
\end{equation*}%
Hence, using (\ref{5.44}), we obtain%
\begin{equation*}
\left\Vert \widetilde{v}_{\varepsilon }^{\ast }-V_{\min ,\lambda
_{4}}\right\Vert _{H_{2}^{1}\left( \Omega \right) \times L_{2,2}\left(
-d,d\right) }\leq C_{3}\delta .
\end{equation*}%
Denote
\begin{equation}
\widehat{V}_{\min ,\lambda _{4}}=V_{\min ,\lambda _{4}}+F.  \label{5.48}
\end{equation}%
We have
\begin{equation*}
\hspace{-1 cm}
\left.
\begin{array}{c}
\left\Vert v_{\varepsilon }^{\ast }-\widehat{V}_{\min ,\lambda
_{4}}\right\Vert _{H_{2}^{1}\left( \Omega \right) \times L_{2,2}\left(
-d,d\right) }= \\
=\left\Vert \left( v_{\varepsilon }^{\ast }-F^{\ast }\right) -\left(
\widehat{V}_{\min ,\lambda _{3}}-F\right) +\left( F^{\ast }-F\right)
\right\Vert _{H_{2}^{1}\left( \Omega \right) \times L_{2,2}\left(
-d,d\right) }= \\
=\left\Vert \left( \widetilde{v}_{\varepsilon }^{\ast }-V_{\min ,\lambda
_{4}}\right) +\left( F^{\ast }-F\right) \right\Vert _{H_{2}^{1}\left( \Omega
\right) \times L_{2,2}\left( -d,d\right) }\leq \\
\leq \left\Vert \widetilde{v}_{\varepsilon }^{\ast }-V_{\min ,\lambda
_{4}}\right\Vert _{H_{2}^{1}\left( \Omega \right) \times L_{2,2}\left(
-d,d\right) }+\left\Vert F-F^{\ast }\right\Vert _{H_{2}^{1}\left( \Omega
\right) \times L_{2,2}\left( -d,d\right) }\leq C_{3}\delta .%
\end{array}
\right.
\end{equation*}
Hence,
\begin{equation}
\left\Vert v_{\varepsilon }^{\ast }-\widehat{V}_{\min ,\lambda
_{4}}\right\Vert _{H_{2}^{1}\left( \Omega \right) \times L_{2,2}\left(
-d,d\right) }\leq C_{3}\delta .  \label{5.49}
\end{equation}%
Using (\ref{5.300}), (\ref{5.49}) and the triangle inequality, we obtain%
\begin{equation*}
\left\Vert \widehat{V}_{\min ,\lambda _{4}}\right\Vert _{H_{2}^{1}\left(
\Omega \right) \times L_{2,2}\left( -d,d\right) }\leq \left\Vert
v_{\varepsilon }^{\ast }\right\Vert _{H_{2}^{1}\left( \Omega \right) \times
L_{2,2}\left( -d,d\right) }+C_{3}\delta =R.
\end{equation*}%
Therefore,
\begin{equation}
\widehat{V}_{\min ,\lambda _{4}}\in \overline{B\left(
g_{1},g_{2},g_{3},g_{4},R\right) }.  \label{5.50}
\end{equation}%
On the other hand, let $v_{\min \lambda _{4}}$ be the minimizer of the
functional $J_{\lambda _{3},\gamma }\left( v\right) $ on the set $\overline{
B\left( g_{1},g_{2},g_{3},g_{4},R\right) },$ which is found in Theorem 5.1,%
\begin{equation}
\min_{\overline{B\left( g_{1},g_{2},g_{3},g_{4},R\right) }}J_{\lambda
_{4},\gamma }\left( v\right) =J_{\lambda _{4},\gamma }\left( v_{\min
,\lambda _{4}}\right) .  \label{5.51}
\end{equation}%
Let $\widetilde{v}_{\min ,\lambda _{4}}=v_{\min \lambda _{4}}-F.$ Then $%
\widetilde{v}_{\min ,\lambda _{4}}\in \overline{B_{0}\left( 2R\right) }$ and
by (\ref{5.420})%
\begin{equation}
J_{\lambda _{4},\gamma }\left( V_{\min ,\lambda _{4}}+F\right) \leq
J_{\lambda _{4},\gamma }\left( \widetilde{v}_{\min ,\lambda _{4}}+F\right)
=J_{\lambda _{4},\gamma }\left( v_{\min ,\lambda _{4}}\right) .  \label{500}
\end{equation}%
However, since by (\ref{5.48}) and (\ref{5.50}) $V_{\min ,\lambda _{4}}+F\in
\overline{B\left( g_{1},g_{2},g_{3},g_{4},R\right) },$ then by (\ref{5.51})
we should have
\begin{equation}
J_{\lambda _{4},\gamma }\left( V_{\min ,\lambda _{3}}+F\right) \geq
J_{\lambda _{4},\gamma }\left( v_{\min ,\lambda _{4}}\right) .  \label{501}
\end{equation}
Since the minimizer is unique, then (\ref{500}) and (\ref{501}) imply that $%
V_{\min ,\lambda _{4}}+F=\widehat{V}_{\min ,\lambda _{4}}=v_{\min ,\lambda
_{4}}.$ Thus, (\ref{5.49}) implies (\ref{5.37}). Estimate (\ref{5.38})
obviously follows from (\ref{5.7}) and (\ref{5.37}). $\square $

\subsection{Global convergence of the gradient descent method}

\label{sec:5.3}

Let $\lambda =\lambda _{4},$ where $\lambda _{4}$ is defined in (\ref{5.35}%
). Consider two sets
\begin{equation*}
\widetilde{B}^{\ast }=B\left( g_{1}^{\ast },g_{2}^{\ast },g_{3}^{\ast
},g_{4}^{\ast },R/3\right) ,\ \widetilde{B}=B\left(
g_{1},g_{2},g_{3},g_{4},R/3\right) .
\end{equation*}%
We assume now that
\begin{equation}
v_{\varepsilon }^{\ast }\in \widetilde{B}^{\ast },v_{\min ,\lambda _{4}}\in
\widetilde{B}.  \label{5.52}
\end{equation}%
Consider the gradient descent method of the minimization of the functional $%
J_{\lambda _{4},\gamma }.$ Consider an arbitrary point
\begin{equation}
v_{0}\in \widetilde{B}.  \label{5.53}
\end{equation}%
Let $\rho >0$ be a small number. Define the sequence of the gradient descent
method as:%
\begin{equation}
v_{m}=v_{m-1}-\rho J_{\lambda _{4},\gamma }^{\prime }\left( v_{m-1}\right)
,m=1,2,...  \label{5.54}
\end{equation}%
Note that since by Theorem 5.1 $J_{\lambda _{4},\gamma }^{\prime }\left(
v_{m-1}\right) \in S_{0}$ for all $m\geq 1,$ then it follows from (\ref{5.5}%
) and (\ref{5.54}) that boundary conditions (\ref{3.10}), (\ref{3.11}) are
kept the same for all vector functions $v_{m}.$ Theorem 5.3 follows
immediately from a combination of Theorem 5.2 with Theorem 6 of \cite{KSAR}.

\textbf{Theorem 5.3}. \emph{Let }$\lambda =\lambda _{4}$\emph{\ and let
conditions of Theorem 5.2, (\ref{5.52}) and (\ref{5.53}) hold. Then there
exists a sufficiently small number }$\rho _{0}\in \left( 0,1\right) $\emph{\
such that for every }$\rho \in \left( 0,\rho _{0}\right) $\emph{\ there
exists a number }$\theta =\theta \left( \rho \right) \in \left( 0,1\right) $%
\emph{\ such that the sequence }$\left\{ v_{m}\right\} _{m=0}^{\infty
}\subset \widetilde{B}$\emph{\ and the following convergence estimates for
the gradient descent method (\ref{5.53}), (\ref{5.54}) hold:}%
\begin{eqnarray*}
& \left\Vert v_{m}-v_{\varepsilon }^{\ast }\right\Vert _{H_{2}^{1}\left(
\Omega \right) \times L_{2,2}\left( -d,d\right) }\leq C_{3}\delta +\theta
^{m}\left\Vert v_{\min ,\lambda _{4}}-v_{0}\right\Vert _{S},\ m=1,...,
\\
&\hspace{0.5 cm} \left\Vert a_{m}-a_{\varepsilon }^{\ast }\right\Vert
_{L_{2}\left( \Omega \right) }\leq C_{3}\delta +\theta ^{m}\left\Vert
v_{\min ,\lambda _{3}}-v_{0}\right\Vert _{S},\ m=1,...,
\end{eqnarray*}%
\emph{where the function }$a_{m}\left( x\right) $\emph{\ is constructed from
the vector function }$v_{m}=\left( p_{m},q_{m}\right) $\emph{\ by the right
hand side of formula (\ref{5.7}), in which }$p_{\varepsilon ,\mbox{comp}}$%
\emph{\ is replaced with }$p_{m}.$

\textbf{Remarks 5.1:}

\begin{enumerate}
\item \emph{Since smallness assumptions are not imposed on the number }$R>0$
\emph{\ and since }$v_{0}\in \widetilde{B}$\emph{\ is an arbitrary point,
then Theorem 5.3 guarantees the global convergence of the gradient descent
method (\ref{5.53}), (\ref{5.54}).}

\item \emph{Even though the requirement of our theory is that the parameter }
$\lambda $\emph{\ of the Carleman Weight Function }$e^{2\lambda z^{2}}$\emph{%
\ \ should be sufficiently large, we have observed in computational
experiments of section 6 that }$\lambda =5$\emph{\ is sufficient, which is
the same as in two previous publications of this group \cite{KTR,RRTE}.
Similar observations about reasonable values of }$\lambda \in \left[ 1,3 %
\right] $\emph{\ were made in other publications about the convexification
method \cite{KL,KSAR,KHJ,KLZ}. On the other hand, we observe in numerical
experiments of section 6 that too large values of }$\lambda =20$ \emph{do
not work well, see Figure \ref{plot_Epsilon0P01_Beta0P001_diffLambda}. This
is because the Carleman Weight Function changes well too rapidly for }$
\lambda =20.$

\item \emph{Conceptually, the considerations of item 2 are similar with
asymptotic theories. Indeed, an asymptotic theory usually claims that if a
certain parameter }$X$\emph{\ is sufficiently large, then a certain formula }
$Y$\emph{\ is valid with a good accuracy. However, in a practical
computation, which always has a specific ranges of parameters, only
numerical experiments can establish reasonable values of }$X,$ \emph{\ for
which }$Y$\emph{\ is valid with a good accuracy. Besides, it is well known
that too large values of }$X$ \emph{often do not work well for numerical
studies. }
\end{enumerate}

\section{Numerical Studies}

\label{sec:6}

We have conducted numerical studies in the 2D case. In our numerical testing
the domain $\Omega $ and the line $\Phi _{d}$ in (\ref{2.1})-(\ref{2.40})
are:
\begin{equation}
\eqalign{
\Omega & =\left\{ x:x_{1}\in \left( -B,B\right) ,z\in \left( a,b\right)
\right\} ,\ B=1/2,\ a=1,\ b=2, \\
& \Phi _{d}=\left\{ x_{\alpha }=(\alpha ,0):\alpha \in \lbrack -d,d]\right\}
,\ d=1/2.
}
\label{6.1}
\end{equation}%
We take $\sigma =0.05$ in the function $f_{\sigma }(x)$ in (\ref{2.5}) and %
(\ref{2.50}). As to the kernel $G(x,\alpha ,\beta )$ in (\ref{2.7}), we
choose the 2-dimensional Henyey-Greenstein function \cite{Heino}:
\begin{equation}
G(x,\alpha ,\beta )=H(\alpha ,\beta )=\frac{1}{2d}\left[ \frac{1-c_{g}^{2}}{
1+c_{g}^{2}-2c_{g}\cos (\alpha -\beta )}\right] .  \label{6.2}
\end{equation}%
Here $c_{g}$ represents the ballistic with $c_{g}=0$ and isotropic
scattering with $c_{g}=1$ \cite{HT1,HT2,HT3}, respectively. In this paper,
we choose $c_{g}=1/2$.

We have chosen the absorption and scattering coefficients $\mu _{a}(x),\mu
_{s}(x)$ in (\ref{2.10})-(\ref{2.12}) as:
\begin{eqnarray}
\mu _{s}(x)& =5,\ x\in \Omega ,\quad \mu _{s}(x)=0,\ x\in \mathbb{R}
^{2}\setminus \Omega .  \label{6.3} \\
\mu _{a}(x)& =\left\{
\begin{array}{cc}
c_{a}=const.>0, & \mbox{inside the tested inclusion,} \\
0, & \mbox{outside the tested inclusion.}%
\end{array}
\right.  \label{6.4}
\end{eqnarray}%
In the numerical tests below, we take $c_{a}=5,10,15,20,30$, and the
inclusions with the shape of the letters `A', `$\Omega $' and `SZ'.

\textbf{Remark 6.1. }\emph{We have intentionally chosen abnormalities with
the shapes of letters to demonstrate that our reconstruction technique works
well for truly hard cases of non-convex abnormalities containing voids.}

By (\ref{2.10}), (\ref{6.3}) and (\ref{6.4}), we have
\begin{equation}
\mbox{inclusion/background contrast}=1+c_{a}/5.  \label{6.5}
\end{equation}%
Following (\ref{6.5}), we define the computed inclusion/background contrast
as:%
\begin{equation}
\hspace{-2 cm}
\mbox{computed inclusion/background contrast}=1+\max (\mbox{computed }\mu
_{a}(x))/5.  \label{6.50}
\end{equation}

\subsection{Data generation}

\label{sec:6.1}

To generate the boundary data (\ref{2.15}) and then (\ref{3.10}), (\ref{3.11}%
) for our CIP, we have solved the Forward Problem posed in section 2. Using
Theorem 2.1, we have solved this problem numerically via the solution of the
linear integral equation (\ref{2.025}) with the condition (\ref{2.026}).

Consider the partition of the domain $\Omega $ and the line $\Phi _{d}$ in %
(\ref{6.1}) with the given mesh sizes $h_{x_{1}},h_{z},h_{\alpha }$:
\begin{eqnarray}
x_{1,i}& =-B+ih_{x_{1}},\quad i=0,1,\cdots ,n_{x_{1}},\quad
n_{x_{1}}=2B/h_{x_{1}},  \label{6.7} \\
z_{j}& =a+jh_{z},\quad j=0,1,\cdots ,n_{z},\quad n_{z}=(b-a)/h_{z},
\label{6.8} \\
\alpha _{k}& =-d+kh_{\alpha },\quad k=0,1,\cdots ,n_{\alpha },\quad
n_{\alpha }=2d/h_{\alpha }.  \label{6.9}
\end{eqnarray}%
Then the set of the discrete points are given as
\begin{eqnarray}
x^{h}& =\left\{ x_{i,j}=(x_{1,i},z_{j}),\ i=0,1,\cdots ,n_{x_{1}},\
j=0,1,\cdots ,n_{z}\right\} ,  \label{6.10} \\
&\hspace{2 cm} \alpha ^{h}=\left\{ \alpha _{k},\ k=0,1,\cdots ,n_{\alpha
}\right\} .  \label{6.11}
\end{eqnarray}%
We have used the grid step sizes $h_{x_{1}}=h_{z}=h_{\alpha }=1/40.$ To
obtain the numerical solution $u(x^{h},\alpha ^{h})$ of the Forward Problem,
we have solved the corresponding linear algebraic system by the Matlab
backslash operator `$\backslash $'. This way we have generated the boundary
data \eref{2.15}. Then, using considerations of subsection 3.1, we have
obtained the boundary data (\ref{3.10}), (\ref{3.11}).

\subsection{Numerical results for the inverse problem}

\label{sec:6.2}

For the inverse problem, we set $h_{x_{1}}=h_{z}=h_{\alpha }=1/20$ to
generate the discrete points in \eref{6.7}-\eref{6.11}. The discrete form
of functional \eref{5.4} is
\begin{equation}
\hspace{-1 cm}
\left.
\begin{array}{c}
J_{\lambda ,\gamma }\left( p^{h},q^{h}\right) = e^{-2\lambda
b^{2}}\int\limits_{-d}^{d}\int\limits_{\Omega }\left\{ \left[
L_{1}^{h}\left( p^{h},q^{h}\right) \right] ^{2}+\left[ L_{2}^{h}\left(
p^{h},q^{h}\right) \right] ^{2}\right\} e^{2\lambda z^{2}}dxd\alpha \\
\hspace{-3.8 cm} + \gamma \left\Vert \left( p^{h},q^{h}\right) \right\Vert
_{S}^{2},%
\end{array}
\right.  \label{6.12}
\end{equation}
where the pair $\left( p^{h},q^{h}\right) (x^{h},\alpha ^{h})$ is the pair
of functions $\left( p,q\right) \left( x,\alpha \right) $ written on the
discrete grid and $L_{1}^{h}\left( p^{h},q^{h}\right) $ and $L_{2}^{h}\left(
p^{h},q^{h}\right) $ are operators $L_{1}$ and $L_{2}$ in (\ref{3.8}) and (%
\ref{3.9}), in which differential operators are written in finite
differences and integrals are written in discrete forms using the
trapezoidal rule.

\begin{figure}[tbph]
\centering
\includegraphics[width = 4in]{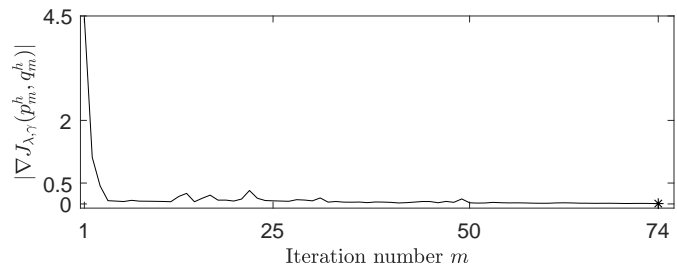}
\caption{Test 1. A typical dependence of $\left\vert \protect\nabla J_{
\protect\lambda ,\protect\gamma }\left( p_{m}^{h},q_{m}^{h}\right)
\right\vert $ from the iteration number $m$. $\ast $ indicates the iteration
number when iterations were stopped because the criterion (\protect\ref{100}%
) was met.}
\label{grad_J}
\end{figure}

To numerically solve the Minimization Problem posed in section 5, we have
minimized functional (\ref{6.12}) with respect to the values of discrete
functions $p_{\varepsilon }^{h}(x^{h},\alpha ^{h})$, $q_{\varepsilon
}^{h}(x^{h},\alpha ^{h})$ at grid points. The Dirichlet boundary conditions
in \eref{3.10} and \eref{3.11} are given as
\begin{equation}
\eqalign{
p_{\varepsilon }^{h}(x_{1,i},z_{j},\alpha )& =g_{1}(x_{1,i},z_{j},\alpha
),\quad q_{\varepsilon }^{h}(x_{1,i},z_{j},\alpha
)=g_{2}(x_{1,i},z_{j},\alpha ), \\
&\hspace{1 cm} i=0,N_{x_{1}},\ j=0,N_{z}.
}
\label{6.13}
\end{equation}%
By the finite difference method, the Neumann boundary conditions in %
\eref{3.10} and \eref{3.11} are given as
\begin{equation}
\eqalign{
-4p_{\varepsilon }^{h}(x,z_{N_{z}-1},\alpha )+p_{\varepsilon
}^{h}(x,z_{N_{z}-2},\alpha )& =2h_{z}g_{3}(x,b,\alpha )-3g_{1}(x,b,\alpha ),
\\
-4q_{\varepsilon }^{h}(x,z_{N_{z}-1},\alpha )+q_{\varepsilon
}^{h}(x,z_{N_{z}-2},\alpha )& =2h_{z}g_{4}(x,b,\alpha )-3g_{2}(x,b,\alpha ).
}
\label{6.14}
\end{equation}%
We have adopted the Matlab's built-in optimization toolbox \textbf{fmincon}
to minimize the function $J_{\lambda ,\gamma }\left( p_{\varepsilon
}^{h},q_{\varepsilon }^{h}\right) $ in \eref{6.12} with the boundary
conditions \eref{6.13} and \eref{6.14}. Here, \eref{6.13} and \eref{6.14}
are the constraint conditions used in each iteration of \textbf{fmincon} to
ensure that the functions $p_{\varepsilon }^{h}(x^{h},\alpha ^{h})$, $%
q_{\varepsilon }^{h}(x^{h},\alpha ^{h})$ at every iteration satisfy the
boundary conditions \eref{3.10} and \eref{3.11}. The iterations of \textbf{%
\ fmincon} were stopped at the iteration number $m$ at which
\begin{equation}
\left\vert \nabla J_{\lambda ,\gamma }\left( p_{m}^{h},q_{m}^{h}\right)
\right\vert <10^{-2},  \label{100}
\end{equation}%
see Figure \ref{grad_J}.

To solve the minimization problem, we need to provide the starting point $%
p_{0}^{h},q_{0}^{h}$ for iterations. With the boundary conditions %
\eref{3.10} and \eref{3.11}, for every $\alpha \in \lbrack -d,d]$, we have
the value of functions $p,q$ on the boundary $\partial \Omega $ via
functions $g_{1},g_{2}$. Using the linear interpolations of boundary
conditions $g_{1},g_{2}$ with respect to $x_{1}-$direction and $z-$%
direction, the initial guess $\left( p_{0},q_{0}\right) $ for the pair of
functions $\left( p,q\right) $ in the domain $\Omega $ is:
\begin{equation}
\left.
\begin{array}{c}
p_{0}(x_{1},z,\alpha ) =\frac{1}{2}\left( \frac{(B-x_{1})}{2B}
g_{1}(-B,z,\alpha )+\frac{(x_{1}+B)}{2B}g_{1}(B,z,\alpha )\right) \\
\hspace{2.0 cm} +\frac{1}{2}\left( \frac{(b-z)}{b-a}g_{1}(x_{1},a,\alpha )+%
\frac{(z-a)}{b-a }g_{1}(x_{1},b,\alpha )\right) , \\
q_{0}(x_{1},z,\alpha ) =\frac{1}{2}\left( \frac{(B-x_{1})}{2B}
g_{2}(-B,z,\alpha )+\frac{(x_{1}+B)}{2B}g_{2}(B,z,\alpha )\right) \\
\hspace{2.0 cm} +\frac{1}{2}\left( \frac{(b-z)}{b-a}g_{2}(x_{1},a,\alpha )+%
\frac{(z-a)}{b-a }g_{2}(x_{1},b,\alpha )\right) .%
\end{array}
\right.  \label{6.15}
\end{equation}
Then the starting point for the minimization of functional (\ref{6.15}) is $%
p_{0}^{h}=p_{0}(x^{h},\alpha ^{h})$, $q_{0}^{h}=q_{0}(x^{h},\alpha ^{h})$.
Even though the first guess does not satisfy required Neumann boundary
conditions in \eref{3.10} and \eref{3.11}, still all follow up iterations
of \textbf{fmincon} satisfy both required boundary conditions: Dirichlet and
Neumann, by the constraint conditions \eref{6.13} and \eref{6.14}.

We introduce the random noise in the boundary data $g(x,\alpha )$ in %
\eref{2.15} on the boundary $\partial \Omega $ as follows:
\begin{equation}
g(x,\alpha )=g(x,\alpha )\left( 1+\delta \zeta _{x}\right) ,  \label{6.20}
\end{equation}%
where $\zeta _{x}$ is the uniformly distributed random variable in the
interval $[0,1]$ depending on the point $x\in \partial \Omega $ with $\delta
=0.03$ and $\delta =0.05$, which correspond respectively to $3\%$ and $5\%$
noise level. Hence, the random noise is also introduced in boundary
conditions $g_{1}\left( x,\alpha \right) ,g_{2}\left( x,\alpha \right)
,g_{3}\left( x,\alpha \right) $ and $g_{4}\left( x,\alpha \right) $ in (\ref%
{3.10}), (\ref{3.11}). These functions are defined via $g\left( x,\alpha
\right) $ using (\ref{3.3}), (\ref{3.6}), (\ref{3.66}) and (\ref{3.67}). We
now explain how did we differentiate the noisy data for $g\left( x,\alpha
\right) $ with respect to $\alpha $ in (\ref{3.6}) and (\ref{3.67}). The
observation data $g(x,\alpha )$ in \eref{2.15} at the boundary $\partial
\Omega $ is generated by the source function $f(x-x_{\alpha })$, whose
position $x_{\alpha }$ is determined by the value $\alpha $ in \eref{2.40}.
Then, for each given $\alpha $, we obtain the corresponding observation data
$g(x,\alpha )$ generated by the source function $f(x-x_{\alpha })$ as well
as the the sample of the random variable $\zeta _{x}$. Since the samples of
the random variable $\zeta _{x}$ for each $\alpha _{k}$ in (\ref{6.9}) are
different, then we use the finite difference method to calculate numerically
the derivative of the noisy data $g(x,\alpha )$ with respect to $\alpha $
with the above mentioned grid step size $h_{\alpha }=1/40.$ Results of the
differentiation were good enough, and we did not observe instabilities.

\textbf{Test 1}. We consider the coefficient $a(x),$ which corresponds to $%
\mu _{a}(x)$ in \eref{6.4} with $c_{a}=5$ inside of the letter `$A$'. The
goal of this test is to find the optimal values of the parameters $%
\varepsilon $ and $\lambda $ for the minimization problem. Noise (\ref{6.20}%
) in the data is not added.

We set $\gamma =0.001$, $\lambda =5$ and perform the numerical tests with
different values of $\varepsilon $. The results are displayed in Figure \ref%
{plot_lambda5_Beta0P001_diffEpsilon}. The reconstruction of the header of
the letter `$A$' is not good for $\varepsilon \in \left[ 0.05,0.2\right] $.
The reconstruction quality improves when $\varepsilon $ decreases while
ranging from 0.2 to 0.01. While $\varepsilon $ varies from 0.01 to 0.001,
the difference between the reconstructions is very small. On the other hand,
when we choose $\varepsilon =0.0001,$ the reconstruction quality of the top
of letter `$A$' becomes worse. In conclusion, although $\varepsilon $ should
be small enough, but not too small. Thus, we choose $\varepsilon =0.01$ as
an optimal value and use this one in all other tests.

\begin{figure}[tbph]
\centering
\includegraphics[width =4in]
	{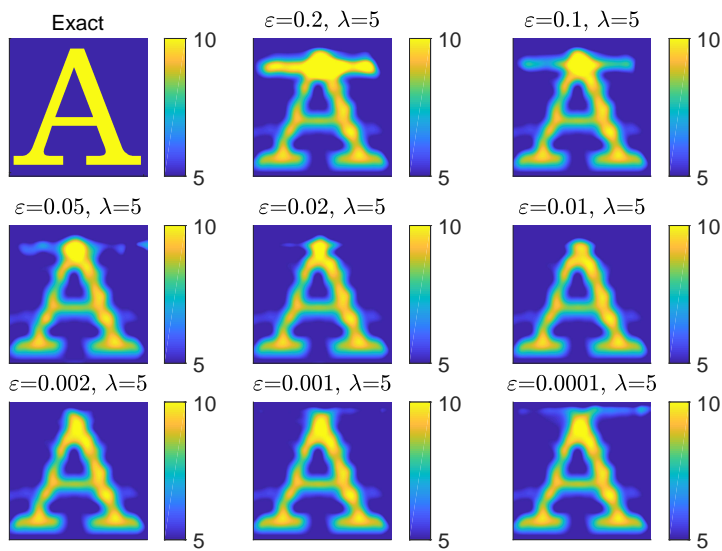}
\caption{Test 1. The reconstructed coefficient $a(x) $, where the function $%
\protect\mu _{a}(x) $ is given in \eref{6.4} with $c_{a}=5$ inside of the
letter `A'. We choose $\protect\lambda=5, \protect\gamma=0.001$, and test
different values of $\protect\varepsilon$. The reconstruction is not good
for $\protect\varepsilon \in \left[ 0.05,0.2\right] $. On the other hand,
the quality of the reconstruction improves when $\protect\varepsilon $
decreases while ranging from 0.2 to 0.01. While $\protect\varepsilon $
varies from 0.01 to 0.001, the difference between the reconstructions is
very small. On the other hand, the choice $\protect\varepsilon =0.0001$
leads to a lower quality reconstruction quality of the top of letter `$A$'.
Thus, we choose $\protect\varepsilon =0.01$ as an optimal value and use this
value in all follow up tests.}
\label{plot_lambda5_Beta0P001_diffEpsilon}
\end{figure}

\begin{figure}[tbph]
\centering
\includegraphics[width =4in]
	{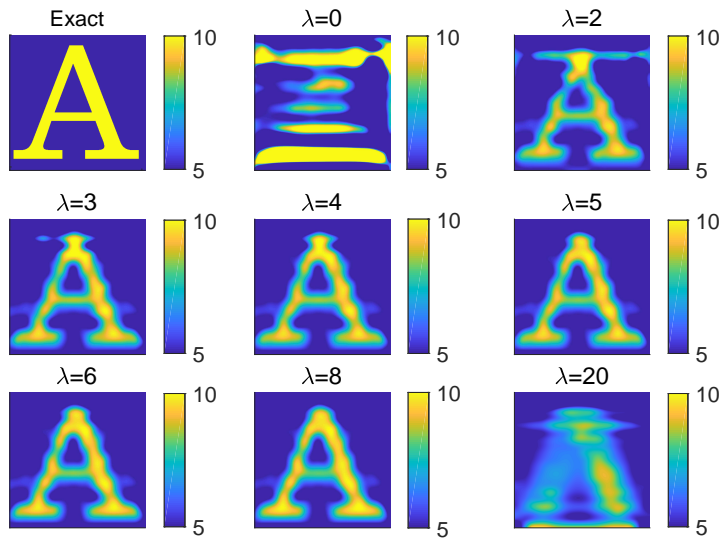}
\caption{Test 1. The reconstructed coefficient $a(x) $, where the function $%
\protect\mu _{a}(x) $ is given in \eref{6.4} with $c_{a}=5$ inside of the
letter `A'. We choose $\protect\gamma=0.001, \protect\varepsilon=0.01$, and
test different values of $\protect\lambda$. The reconstructions are
unsatisfactory when $\protect\lambda $ is too small, $\protect\lambda=0, 2$.
Next, the reconstruction quality is improved when $\protect\lambda$ changes
between 3 and 5, and it is stabilized for $\protect\lambda =5,6,8$. On the
other hand, the reconstruction quality deteriorates at $\protect\lambda =20$%
. Thus, we choose $\protect\lambda =5$ as the optimal value. }
\label{plot_Epsilon0P01_Beta0P001_diffLambda}
\end{figure}

Now we want to select an optimal value of the parameter $\lambda .$ We take $%
\gamma =0.001,\varepsilon =0.01$ and test values $\lambda =0,2,3,4,5,6,8,20 $%
. Results are presented on Figure \ref{plot_Epsilon0P01_Beta0P001_diffLambda}%
. The parameter $\lambda $ cannot be neither too small nor too large. The
reconstructions are unsatisfactory for $\lambda =0,2,20$. On the other hand,
the reconstructions become better when $\lambda $ ranges from 3 to 5, and
they are stabilized for $\lambda =5,6,8$. Thus, we choose $\lambda =5$ as
the optimal value, see items 2 and 3 of Remarks 5.1 for a relevant
discussion.

In summary, we use in the tests below
\begin{equation}
\gamma =0.001,\varepsilon =0.01,\lambda =5.  \label{101}
\end{equation}

\begin{figure}[tbph]
\centering
\includegraphics[width = 5in]{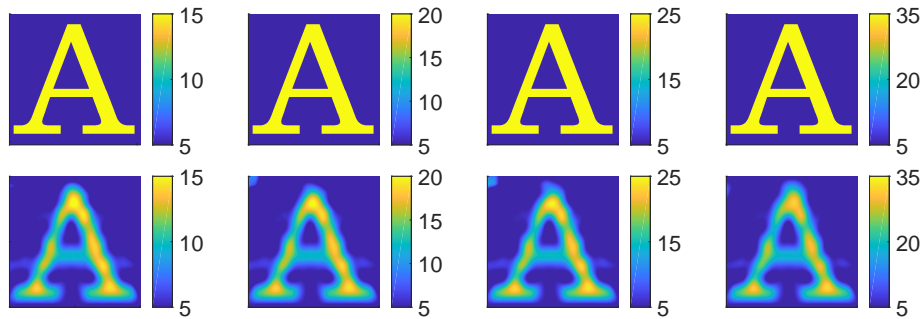}
\caption{Test 2. Exact (top) and reconstructed (bottom) coefficient $a(x)$
with $c_{a}=10,15,20,30$ inside of the letter `A' as in (\protect\ref{6.4})
for the values of the parameters being the same as in (\protect\ref{101}).
The inclusion/background contrasts in (\protect\ref{6.50}) are respectively $%
3:1$ , $4:1$, $5:1$ and $7:1$. Even though there is a small deterioration at
$c_{a}=20$ and $c_{a}=30$, the accuracy of the reconstruction remains
basically the same for these four choices of $c_{a}$. The computed
inclusion/background contrasts in (\protect\ref{6.50}) are accurate..}
\label{plot_A10_A15_A20_A30}
\end{figure}

\textbf{Test 2}. We take the same values of parameters as listed in (\ref%
{101}). We consider the coefficient $a(x)$ corresponding to $\mu _{a}(x)$ in %
\eref{6.4} with $c_{a}=10,15,20,30$ inside of the letter `$A$'. Hence, the
inclusion/background contrasts in (\ref{6.5}) are respectively $3:1$, $4:1$,
$5:1$ and $7:1$. Noise (\ref{6.20}) in the data is not added. The results
are displayed in Figure \ref{plot_A10_A15_A20_A30}. The reconstruction
quality is good for these four cases, although it slightly deteriorates at $%
c=20$ and $c=30$. The computed inclusion/background contrasts (\ref{6.50})
are accurate.

\textbf{Test 3}. We use the same values of parameters as ones in (\ref{101}%
). We consider the coefficient $a(x)$ corresponding to $\mu _{a}(x)$ in %
\eref{6.4} with $c_{a}=5$ inside of the letter `$\Omega $'. Noise (\ref%
{6.20}) in the data is not added. The result is displayed in Figure \ref%
{plot_Omega}. The reconstruction is quite accurate.

\begin{figure}[tbph]
\centering
\includegraphics[width = 3.5in]{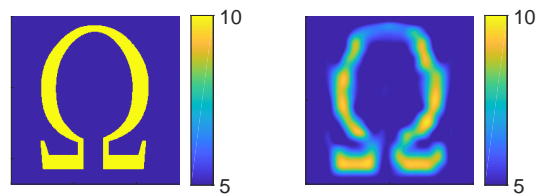}
\caption{Test 3. Exact (left) and reconstructed (right) coefficient $a(x)$
with the shape of the letter `$\Omega $' with $c_{a} = 5$ in \eref{6.4}.
The values of the parameters are the same as ones in (\protect\ref{101}).
The reconstruction is accurate, also, see (\protect\ref{6.5}) and (\protect
\ref{6.50}).}
\label{plot_Omega}
\end{figure}

\begin{figure}[tbph]
\centering
\includegraphics[width = 3.5in]{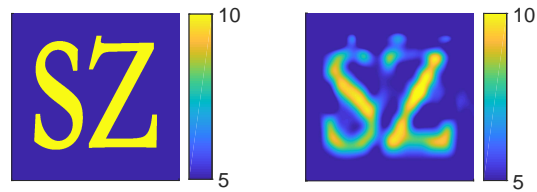}
\caption{Test 4. Exact (left) and reconstructed (right) coefficient $a(x)$
for $\protect\mu_{a}(x)$ in \eref{6.4} with $c_{a} = 5$ inside of two
letters `SZ'. The values of the parameters are the same as ones in (\protect
\ref{101}).The reconstruction is worse than the one for the case of the
single letter `$\Omega $' in Figure \protect\ref{plot_Omega}. Nevertheless,
the reconstruction is still good and the computed inclusion/background
contrasts in (\protect\ref{6.50}) are accurate in both letters.}
\label{plot_SZ}
\end{figure}

\textbf{Test 4}. We consider the coefficient $a(x)$ corresponding to $\mu
_{a}(x)$ in \eref{6.4} with $c_{a}=5$ inside of two letters `SZ', which are
two letters in the name of the city (Shenzhen) were the second and the third
authors reside. Noise (\ref{6.20}) in the data is not added. Results are
exhibited in Figure \ref{plot_SZ}. The reconstruction is worse than the one
for the case of the single letter `$\Omega $' in Figure \ref{plot_Omega}.
Nevertheless, the reconstruction is still good and the computed
inclusion/background contrasts in (\ref{6.50}) are both accurate in these
two letters.

\textbf{Test 5}. We now consider the noisy data, as in \eref{6.20}, with $%
\delta =0.03$ and $\delta =0.05,$ i.e. with 3\% and 5\% noise level. We
reconstruct the coefficient $a(x)$ with the shape of the letters `A' and `$%
\Omega $' corresponding to $\mu _{a}(x)$ in \eref{6.4} with $c_{a}=5$
inside of two letters. The results are displayed in Figure \ref%
{plot_AddNoise}. In all these four cases, reconstructions of shapes of
inclusions and the inclusion/background contrasts in (\ref{6.50}) are
accurate.

\begin{figure}[tbph]
\centering
\includegraphics[width = 5in]{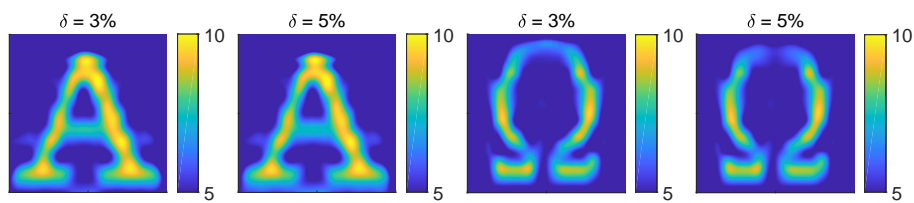}
\caption{Test 5. Reconstructed coefficient $a(x)$ with the shape of letters
`A' and `$\Omega $' with $c_{a}=5$ from noisy data (\protect\ref{6.20}) with
$\protect\delta =0.03$ and $\protect\delta =0.05$, i.e. with 3\% and 5\%
noise level. The values of the parameters are the same as ones in (\protect
\ref{101}). In all these four cases, both reconstructions and the
inclusion/background contrasts in (\protect\ref{6.50}) are accurate.}
\label{plot_AddNoise}
\end{figure}

\section*{References}


\end{document}